\documentstyle{article}

\def\pr{\mbox{\small \O}}
\def\pl{\!+\!}
\def\mn{\!-\!}
\def\prop#1#2{\noindent{\bf #1.} {#2}}

\def\dkz{\vspace{2ex} \noindent{\em Proof. }}
\def\qed{\hfill $\dashv$ \vspace{2ex}}

\begin{document}

\title{Shuffles and Concatenations\\ in Constructing of Graphs}
\author{{\sc Kosta Do\v sen} and {\sc Zoran Petri\' c}
\\[1ex]
{\small Mathematical Institute, SANU}\\[-.5ex]
{\small Knez Mihailova 36, p.f.\ 367, 11001 Belgrade,
Serbia}\\[-.5ex]
{\small email: \{kosta, zpetric\}@mi.sanu.ac.rs}}
\date{}
\maketitle

\begin{abstract}
\noindent This is an investigation of the role of shuffling and
concatenating in the theory of graph drawing. A simple syntactic
description of these and related operations is proved complete in
the context of finite partial orders, as general as possible. An
explanation based on that is given for a previously investigated
collapse of the permutohedron into the associahedron, and for
collapses into other less familiar polyhedra, including the
cyclohedron. Such polyhedra have been considered recently in
connection with the notion of tubing, which is closely related to
tree-like finite partial orders defined simply and investigated
here in detail. Like the associahedron, some of these other
polyhedra are involved in categorial coherence questions, which
will be treated in a sequel to this paper.
\end{abstract}


\vspace{.5ex}


\section{Introduction}
Shuffles and concatenations, which are usually considered only for
finite linear orders, are here defined for arbitrary binary
relations (see Section~4). Shuffles serve to define on sets of
relations an associative and commutative partial operation, which
we call shuffle sum; concatenations serve analogously to define on
sets of relations an associative partial operation, which we call
concatenation product.

Shuffle sum and concatenation product are interesting to us for
the following reason. The one-one map $L$, which assigns to a
partial order all its linear extensions, maps disjoint union and
concatenation of partial orders into shuffle sum and concatenation
product respectively (see Section~4). And here is why disjoint
union and concatenation of partial orders are interesting to us.

We associate with a given graph $\Gamma$ a set of terms
representing tree-like finite partial orders $T(\Gamma)$, each of
which may be understood as a possible history of the constructing,
or, in reverse order, destructing, of $\Gamma$. The set
$T(\Gamma)$ determines the graph $\Gamma$ uniquely, i.e.\ the map
$T$ is one-one. The tree-like partial orders of $T(\Gamma)$ are
closely related to the \emph{tubings} of \cite{CD06}, but they are
defined more simply.

The members of $T(\Gamma)$ are built inductively in a simple
manner with the help of two operations corresponding to disjoint
union and concatenation. These operations correspond via the map
$L$ mentioned above to shuffle sum and concatenation product.

The members of $T(\Gamma)$ label vertices of polyhedra that are
obtained from permutohedra by collapsing connected families of
vertices into a single vertex. We use the map $L$ to assign to a
member of $T(\Gamma)$ the permutations in a connected family of
vertices of the permutohedron, which are collapsed into a single
vertex. The collapsing in question that produces associahedra has
been studied previously in \cite{T97}. For a suitable choice of
$\Gamma$, we obtain a collapsing that produces cyclohedra, and
other choices yield less familiar polyhedra.

These polyhedra stand for commuting diagrams that arise in various
coherence questions in category theory. It is shown in \cite{DP06}
how Mac Lane's pentagon of monoidal coherence arises by a
collapsing of the same kind we have here from a hexagon involved
in symmetric monoidal coherence, and this matter is related to the
collapsing investigated in \cite{T97}.

Some similar coherence questions based on the conceptual apparatus
introduced in \cite{DP06}, which we intend to treat in the future,
involve some of the less familiar polyhedra that occur as examples
in the present paper. The \emph{hemiassociahedron} of Example 5.14
arises in the definition of a coherent notion of weak Cat-operad.
A Cat-operad is an operad enriched over the category Cat of all
small categories, as a 2-category with small hom-categories is a
category enriched over Cat (for the notion of operad see
\cite{L97}). The notion of weak Cat-operad is to the notion of
Cat-operad what the notion of bicategory is to the notion of
2-category. Our notion of weak Cat-operad is coherent in the sense
that all the diagrams of canonical arrows commute, as in Mac
Lane's notion of monoidal category. The commuting diagrams assumed
for this notion may be pasted to make the hemiassociahedron,
besides making the three-dimensional associahedron and
permutohedron. We demonstrate that in \cite{DP10}, for which the
present paper lays the ground.

Our examples of collapsing depend on specific graphs $\Gamma$, but
we show that we have a general phenomenon, not to be found only in
our examples. The maps $T$ and $L$ for a given graph $\Gamma$ with
$n$ vertices induce an equivalence relation on the set of vertices
of the $n\mn 1$-dimensional permutohedron (see Section~5).

Our tree-like partial orders are easily described syntactically
with two partial binary operations, one, corresponding to disjoint
union, associative and commutative, and the other, corresponding
to concatenation, just associative. This syntactic description
covers a wider class of finite partial orders, with a property
more general than difunctionality, which we call trifunctionality.
We obtain an isomorphism result concerning this matter (see
Sections 2 and 3).

Our treatment of shuffles and concatenations may be connected to
the algebras studied in \cite{JR79} and elsewhere. The connection
is however not clear.

After completing this paper, we learned that shuffles and
concatenations have been studied previously in a manner related to
ours in a number of papers cited in \cite{AF06} (Section 3). In
particular, the results of our Sections 2 and 3 were anticipated
in \cite{G81} and  \cite{VTL82} (Theorem~1; see also \cite{T94},
Lemma~1, and \cite{G88}, Theorem 3.1). Further references to
previous related work may be found in \cite{DP10} and
\cite{DP10a}, which develop the research we have started here, and
give motivation for it.

\section{Disjoint union and concatenation of relations}
In this section we study preliminary matters concerning the
partial operations of disjoint union and concatenation of binary
relations. These operations are partial because we require
disjointness of domains. We are interested in particular in
applying these operations to partial orders that satisfy a
property we call trifunctionality, which generalizes
difunctionality (see references below). The results of this
section prepare the ground for the isomorphism result of the next
section.

A \emph{relation} on a set $X$ is, as usual, an ordered pair
$\langle R,X\rangle$ such that $R\subseteq X^2$. (We deal only
with binary relations in this paper.) The set $X$ is the
\emph{domain} of $\langle R,X\rangle$.

For the relations $\langle R,X\rangle$ and $\langle S,Y\rangle$
such that ${X\cap Y=\pr}$ we have
\begin{tabbing}
\hspace{3em}\=$\langle R,X\rangle$\=$\,+\,$\=$\langle S,Y\rangle
$\=$\;=_{df}\langle R\cup S,X\cup Y\rangle$,\\*[.5ex]

\>$\langle R,X\rangle$\>\hspace{.4em}$\cdot$\>$\langle
S,Y\rangle$\>$\;=_{df}\langle R\cup S\cup (X\times Y),X\cup
Y\rangle$.
\end{tabbing}
The operation $+$ is disjoint union, while $\cdot$ could be called
\emph{concatenation}, because this is what it is when $\langle
R,X\rangle$ and $\langle S,Y\rangle$ are linear orders on finite
domains, i.e.\ finite sequences. It is clear that $+$ is
associative and commutative, while $\cdot$ is associative without
being commutative for $X$ and $Y$ nonempty (for $X$ or $Y$ empty,
$+$ and $\cdot$ coincide). It is easy to verify the following.

\vspace{2ex}

\prop{Remark $+$}{If $\langle R,X\rangle$ is such that $X=X_1\cup
X_2$, $X_1\cap X_2=\pr$ and for every $x_1$ in $X_1$ and every
$x_2$ in $X_2$ we have $(x_1,x_2)\notin R$ and $(x_2,x_1)\notin
R$, then there are relations $\langle R_1,X_1\rangle$ and $\langle
R_2,X_2\rangle$ such that $\langle R,X\rangle = \langle
R_1,X_1\rangle + \langle R_2,X_2\rangle$.}

\vspace{2ex}

\prop{Remark $\,\cdot\;$}{If $\langle R,X\rangle$ is such that
$X=X_1\cup X_2$, $X_1\cap X_2=\pr$ and for every $x_1$ in $X_1$
and every $x_2$ in $X_2$ we have $(x_1,x_2)\in R$ and
$(x_2,x_1)\notin R$, then there are relations $\langle
R_1,X_1\rangle$ and $\langle R_2,X_2\rangle$ such that $\langle
R,X\rangle = \langle R_1,X_1\rangle\cdot\langle R_2,X_2\rangle$.}

\vspace{2ex}

\emph{Partial orders} in this paper will be \emph{strict} partial
orders--i.e.\ relations that are irreflexive and transitive. Note
that if $\langle R,X\rangle$ is a partial order, then in Remark
$\cdot$ we may omit the conjunct $(x_2,x_1)\notin R$, which
follows from $(x_1,x_2)\in R$. We can trivially prove the
following.

\vspace{2ex}

\prop{Proposition 2.1}{If $X\cap Y=\pr$, then $\langle R,X\rangle$
and $\langle S,Y\rangle$ are partial orders iff $\langle
R,X\rangle + \langle S,Y\rangle$ is a partial order, and the same
with $\cdot$ instead of~$+$.}

\vspace{2ex}

We call a relation $\langle R,X\rangle$ \emph{trifunctional} when
for every $x$, $y$, $z$ and $u$ in $X$ we have that if $(x,z)$ and
$(y,z)$ and $(y,u)$ are in $R$, then either $(x,u)$ or $(y,x)$ or
$(u,z)$ is in $R$. The following picture helps to grasp this
implication:
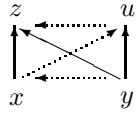
\begin{center}
\begin{picture}(45,49)
\put(2,16){\vector(0,1){20}}

\put(44,16){\vector(0,1){20}}

\put(42,16){\vector(-2,1){38}}

\put(0,42){\makebox(0,0)[l]{$z$}}

\put(42,42){\makebox(0,0)[l]{$u$}}

\put(0,8){\makebox(0,0)[l]{$x$}}

\put(42,8){\makebox(0,0)[l]{$y$}}

\multiput(14,16)(2,0){12}{\line(1,0){.5}}

\multiput(14,36)(2,0){12}{\line(1,0){.5}}

\multiput(5,16.5)(2,1){17}{\line(1,0){.5}}

\put(12,16){\vector(-1,0){2}}

\put(12,36){\vector(-1,0){2}}

\put(40.5,34.5){\vector(3,2){.5}}

\end{picture}
\end{center}
If in this implication we omit the disjuncts $(y,x)\in R$ and
$(u,z)\in R$ from the consequent, then we obtain the implication
that defines difunctional relations (see \cite{R48} and
\cite{SS93}, Section 4.4; our term ``trifunctional'' is motivated
by ``difunctional'', and by the fact that we have three conjuncts
in the antecedent and three disjuncts in the consequent). We can
prove the following.

\vspace{2ex}

\prop{Proposition 2.2}{If $X\cap Y=\pr$, then $\langle R,X\rangle$
and $\langle S,Y\rangle$ are trifunctional relations iff $\langle
R,X\rangle + \langle S,Y\rangle$ is trifunctional, and the same
with $\cdot$ instead of~$+$.}

\dkz For $+$ the proof is trivial, and for $\cdot$ the direction
from right to left is trivial. It remains to prove that if
$\langle R,X\rangle$ and $\langle S,Y\rangle$ are trifunctional,
then $\langle R,X\rangle \cdot \langle S,Y\rangle$ is
trifunctional. So suppose that $\langle R,X\rangle$ and $\langle
S,Y\rangle$ are trifunctional, and suppose that for $x$, $y$, $z$
and $u$ in $X\cup Y$ we have $(x,z)$, $(y,z)$ and $(y,u)$ in
$R\cup S\cup (X\times Y)$. We have the following cases:
\begin{tabbing}
\hspace{1em}1) $z\in X$; then $x,y\in X$, and we have the subcases:\\
\hspace{1.5em}1.1) $u\in X$; then we appeal to the
trifunctionality of $\langle R,X\rangle$;\\
\hspace{1.5em}1.2) $u\in Y$; then we have $(x,u)\in X\times
Y$;\\[1ex]
\hspace{1em}2) $z\in Y$; then we have the subcases:\\
\hspace{1.5em}2.1) $u\in X$; then $(u,z)\in X\times Y$;\\
\hspace{1.5em}2.2) $u\in Y$; then we have the subcases:\\
\hspace{2em}2.21) $x\in X$; then $(x,u)\in X\times Y$;\\
\hspace{2em}2.22) $x\in Y$; then we have the subcases:\\
\hspace{2.5em}2.221) $y\in X$; then $(y,x)\in X\times Y$;\\
\hspace{2.5em}2.222) $y\in Y$; then we appeal to the
trifunctionality of $\langle S,Y\rangle$.\`$\dashv$
\end{tabbing}

For a relation $\langle R,X\rangle$ and $x_1$ and $x_n$, where
$n\geq 2$, distinct elements of $X$ we write $x_1\sim_R x_n$ when
there is a sequence $x_1\ldots x_n$ such that for every
$i\in\{1,\ldots ,n\mn 1\}$ we have $(x_i,x_{i+1})\in R$ or
$(x_{i+1},x_i)\in R$. We say that $\langle R,X\rangle$ is
\emph{connected} when for every two distinct $x$ and $y$ in $X$ we
have $x\sim_R y$. It is trivial to prove the following
proposition.

\vspace{2ex}

\prop{Proposition 2.3}{If for the relations $\langle R,X\rangle$
and $\langle S,Y\rangle$ we have that $X\cap Y=\pr$, $X\neq\pr$
and $Y\neq\pr$, then $\langle R,X\rangle + \langle S,Y\rangle$ is
not connected and $\langle R,X\rangle \cdot \langle S,Y\rangle$ is
connected.}

\vspace{2ex}

By relying on this proposition we obtain easily the following
proposition, which will be applied in the proof of the
completeness Proposition 3.1 in the next section.

\vspace{2ex}

\prop{Proposition 2.4}{If $\langle R_1,X_1\rangle + \ldots +
\langle R_n,X_n\rangle = \langle S_1,Y_1\rangle + \ldots + \langle
S_m,Y_m\rangle$, and for every $i\in\{1,\ldots,n\}$ and every
$j\in\{1,\ldots,m\}$ we have that $R_i$ and $S_j$ are connected,
while $X_i$ and $Y_j$ are not empty, then $n=m$ and there is a
bijection $\pi\!:\{1,\ldots,n\}\rightarrow \{1,\ldots,m\}$ such
that for every $i\in\{1,\ldots,n\}$ we have $\langle
R_i,X_i\rangle = \langle S_{\pi(i)},Y_{\pi(i)}\rangle$.}

\vspace{2ex}

We also have the following proposition.

\vspace{2ex}

\prop{Proposition 2.5}{If for the relations $\langle
R_1,X\rangle$, $\langle R_2,X\rangle$ and $\langle S,Y\rangle$ we
have that $X\cap Y=\pr$, then $\langle R_1,X\rangle \cdot \langle
S,Y\rangle = \langle R_2,X\rangle \cdot \langle S,Y\rangle$ or
$\langle S,Y\rangle\cdot\langle R_1,X\rangle =  \langle
S,Y\rangle\cdot\langle R_2,X\rangle$ implies $R_1=R_2$.}

\dkz Suppose $\langle R_1,X\rangle \cdot \langle S,Y\rangle =
\langle R_2,X\rangle \cdot \langle S,Y\rangle$. If $(x,y)\in R_1$,
then $(x,y)\in R_1\cup S\cup (X\times Y)$, and hence $(x,y)\in
R_2\cup S\cup (X\times Y)$; but then $(x,y)\in R_2$, because
$x,y\in X$. So $R_1\subseteq R_2$, and we demonstrate in the same
manner $R_2\subseteq R_1$. That $\langle S,Y\rangle\cdot\langle
R_1,X\rangle =  \langle S,Y\rangle\cdot\langle R_2,X\rangle$
implies $R_1=R_2$ is demonstrated analogously. \qed

We use this proposition to establish the following proposition,
which will be applied in the proof of the completeness Proposition
3.1 in the next section.

\vspace{2ex}

\prop{Proposition 2.6}{If $\langle R_1,X_1\rangle\cdot\langle
S_1,Y_1\rangle  = \langle R_2,X_2\rangle\cdot\langle
S_2,Y_2\rangle$, where $\langle S_1,Y_1\rangle$ and $\langle
S_2,Y_2\rangle$ are either not connected or their domains are
singletons, then $\langle R_1,X_1\rangle = \langle R_2,X_2\rangle$
and $\langle S_1,Y_1\rangle  = \langle S_2,Y_2\rangle$.}

\dkz Let $\langle T,Z\rangle =\langle R_1,X_1\rangle\cdot\langle
S_1,Y_1\rangle  = \langle R_2,X_2\rangle\cdot\langle
S_2,Y_2\rangle$. It is clear from the definition of $\cdot$ that
we have
\begin{tabbing}
\hspace{3em}\=$(1)$\hspace{1em}\=$(x\in Y_1\;\;$\=$\&\;
(x,y)\in T)\Rightarrow y\in Y_1$,\\[.5ex]
\>$(2)$\>$(x\in X_2$\>$\&\; y\in Y_2)\Rightarrow (x,y)\in T$.
\end{tabbing}
Note also that from the assumption that $\langle S_1,Y_1\rangle$
and $\langle S_2,Y_2\rangle$ are either not connected or their
domains are singletons it follows that $Y_1$ and $Y_2$ are not
empty.

We show by \emph{reductio ad absurdum} that $Y_1\subseteq Y_2$.
Suppose that not $Y_1\subseteq Y_2$; so there is an $x$ in $Y_1$
such that $x\notin Y_2$, which implies that $x\in X_2$. Then we
show that $Y_2\subseteq Y_1$:
\begin{tabbing}
\hspace{3em}$y\in Y_2\;$\=$\Rightarrow (x,y)\in T$, by
(2),\\*[.5ex] \>$\Rightarrow y\in Y_1$, by (1).
\end{tabbing}
The set $Y_1$ cannot be a singleton, because if it were that, then
$Y_2$, which is not empty, would be the same singleton, and we
supposed that we do not have $Y_1\subseteq Y_2$. So $Y_1$ is not a
singleton, and hence $\langle S_1,Y_1\rangle$ is not connected.

Let $y_1$ and $y_2$ be two distinct elements of $Y_1$ such that we
do not have $y_1\sim_{S_1} y_2$. The following three cases exhaust
all the possibilities for $y_1$ and $y_2$ as elements of $X_2\cup
Y_2$.

The first case is when one of $y_1$ and $y_2$ is in $X_2$ and the
other is in $Y_2$. Let $y_1$ be in $X_2$ and $y_2$ in $Y_2$. Then
by (2) we obtain $(y_1,y_2)\in T$, and since $y_1$ and $y_2$ are
in $Y_1$, we have $(y_1,y_2)\in S_1$, which is a contradiction.

The second case is when $y_1$ and $y_2$ are both in $X_2$. Since
$Y_2$ is not empty, for some $y$ in $Y_2$ we have that $(y_1,y)$
and $(y_2,y)$ are in $T$. Since $Y_2\subseteq Y_1$, we have $y\in
Y_1$, from which we infer that $(y_1,y)$ and $(y_2,y)$ are in
$S_1$; this is a contradiction.

The third case is when $y_1$ and $y_2$ are both in $Y_2$. Since
$x\in X_2$, we have that $(x,y_1)$ and $(x,y_2)$ are in $T$, and
since $x\in Y_1$, we have that $(x,y_1)$ and $(x,y_2)$ are in
$S_1$, which is a contradiction.

So we have established that $Y_1\subseteq Y_2$, and we establish
in an analogous manner that $Y_2\subseteq Y_1$. Hence we have
$Y_1=Y_2$.

We show now that $S_1=S_2$. We have first that
\begin{tabbing}
\hspace{3em}$(x,y)\in S_1\;$\=$\Rightarrow (x,y)\in T \;\;\&\;\;
x,y\in Y_1$,\\[.5ex]
\>$\Rightarrow (x,y)\in T \;\;\&\;\;
x,y\in Y_2$,\\[.5ex]
\>$\Rightarrow (x,y)\in S_2$.
\end{tabbing}
So $S_1\subseteq S_2$, and we show analogously that $S_2\subseteq
S_1$.

Let $Y=Y_1=Y_2$. Since $X_1$ and $Y$ are disjoint, and $X_2$ and
$Y$ are disjoint too, from $X_1\cup Y=X_2\cup Y$ we infer
$X_1=X_2$, and then by Proposition 2.5 we conclude that $\langle
R_1,X_1\rangle = \langle R_2,X_2\rangle$. \qed

\section{Diversified \textbf{S}-terms and relations in FTP}
In this section we characterize syntactically in a very simple
manner trifunctional partial orders on finite sets. This is a
freely generated structure, i.e.\ algebra, with two partial
operations, one associative and commutative, corresponding to
disjoint union, and the other associative, corresponding to
concatenation. The operations are partial because we require that
every free generator occurs just once in an element of our
structure. We prove that this syntactically defined structure is
isomorphic to the structure of trifunctional partial orders on
finite sets with the operations of disjoint union and
concatenation.

Consider terms built out of an infinite set of variables, which we
denote by $x$, $y$, $z,\ldots$, $x_1, \ldots$ with the binary
operations $+$ and $\cdot\;$, which we call \emph{sum} and
\emph{product}. Consider structures, i.e.\ algebras, with two
binary operations $+$ and $\cdot$ such that $+$ is associative and
commutative, while $\cdot$ is associative. Let \textbf{S} be the
structure of this kind freely generated by infinitely many
generators. We may take that the elements of \textbf{S} are
equivalence classes of the terms introduced above, which hence we
call \textbf{S}\emph{-terms}, while the variables $x$, $y$,
$z,\ldots$ are \textbf{S}\emph{-variables}. On these equivalence
classes we define the operations $+$ and $\cdot$ by
$[t]+[s]=[t+s]$ and $[t]\cdot[s]=[t\cdot s]$.

An \textbf{S}-term is called \emph{diversified} when no
\textbf{S}-variable occurs in it more than once. Since
associativity and commutativity preserve diversification, it is
clear that if the equivalence class $[t]$ is an element of
\textbf{S} for $t$ a diversified \textbf{S}-term, then every
element of $[t]$ is diversified. We say that the element $[t]$ of
\textbf{S} is \emph{diversified} when $t$ is a diversified
\textbf{S}-term.

Let FTP be the set of trifunctional partial orders on nonempty
finite sets of \textbf{S}-variables. We define by induction on
complexity a map $\kappa$ from the set of diversified
\textbf{S}-terms to the set FTP:
\begin{tabbing}
\hspace{3em}\=$\kappa(x)=\langle \pr,\{x\}\rangle$,\\[.5ex]
\>$\kappa(t + s)\,$\=$=\kappa(t)\,$\=$+\, \kappa(s)$,\\[.5ex]
\>$\kappa(t \cdot s)$\>$=\kappa(t)$\>$\,\cdot\; \kappa(s)$.
\end{tabbing}
That $\kappa(t)$ is indeed a member of FTP for every diversified
\textbf{S}-term $t$ follows from the fact that the relation
$\langle \pr,\{x\}\rangle$ is in FTP, and from Propositions 2.1
and 2.2.

Since the operation $+$ on relations is associative and
commutative, while $\cdot$ is associative, the map $\kappa$
induces a map $K$ from the set of diversified elements of
\textbf{S} to FTP, which is defined by:
\begin{tabbing}
\hspace{3em}$K[t]=\kappa(t)$.
\end{tabbing}
We use $K[t]$ as an abbreviation for $K([t])$. We can prove the
following completeness proposition.

\vspace{2ex}

\prop{Proposition 3.1}{The map $K$ is one-one.}

\dkz Suppose $\kappa(t)=\kappa(s)$. We proceed by induction on the
number $k$ of \textbf{S}-variables in $t$. Since the domains of
the relations $\kappa(t)$ and $\kappa(s)$ are the same, the same
\textbf{S}-variables occur in $t$ and $s$, and hence $k$ is also
the number of \textbf{S}-variables in~$s$.

If $k=1$, then $t$ and $s$ are the same \textbf{S}-variable. If
$k>1$, let $t$ be of the form $t_1+\ldots +t_n$ and $s$ of the
form $s_1+\ldots +s_m$, for $n\geq 1$ and $m\geq 1$, with $t_i$,
for $i\in\{1,\ldots ,n\}$, and $s_j$, for $j\in\{1,\ldots ,m\}$,
\textbf{S}-variables or products. (Since $k>1$, it is impossible
that $n=1$ and $t_1$ is an \textbf{S}-variable.) Since we have
\begin{tabbing}
\hspace{3em}$\kappa(t)=\kappa(t_1)+\ldots
+\kappa(t_n)=\kappa(s_1)+\ldots +\kappa(s_m)=\kappa(s)$,
\end{tabbing}
by Proposition 2.4 we conclude that $n=m$, and that there is a
bijection $\pi\!:\{1,\ldots,n\}\rightarrow \{1,\ldots,m\}$ such
that for every $i\in\{1,\ldots,n\}$ we have
$\kappa(t_i)=\kappa(s_{\pi(i)})$.

If $n=m>1$, then by the induction hypothesis we have
$[t_i]=[s_{\pi(i)}]$ for every $i\in\{1,\ldots,n\}$, and hence
$[t]=[s]$, by the associativity and commutativity of $+$
in~\textbf{S}.

If $n=m=1$, then $t$ and $s$ are products, and by the
associativity of $\cdot$ in \textbf{S} we have $[t]=[t_1\cdot
t_2]$ and $[s]=[s_1\cdot s_2]$ for $t_2$ and $s_2$ either sums or
\textbf{S}-variables. Since we have
\begin{tabbing}
\hspace{3em}$\kappa(t)=\kappa(t_1)\cdot\kappa(t_2)=
\kappa(s_1)\cdot\kappa(s_2)=\kappa(s)$,
\end{tabbing}
by Proposition 2.6 and the induction hypothesis we obtain
$[t]=[s]$. \qed

For the proof of Proposition 3.2 below, which will help us to
establish that $K$, besides being one-one, is also onto, we need
the notion of \emph{inner element} of $X$ for a relation $\langle
R,X\rangle$; this is an element $y$ of $X$ such that for some $x$
and $z$ in $X$ we have $(x,y)\in R$ and $(y,z)\in R$.

For a relation $\langle R,X\rangle$ and $y$ an element of $X$, let
the relation $\langle R\mn y,X\mn\{y\}\rangle$ be defined by
\begin{tabbing}
\hspace{3em}$R\mn y=\{(u,v)\in R\mid u\neq y \;\;\&\;\; v\neq
y\}$.
\end{tabbing}
Then we can formulate the following, which is easy to establish.

\vspace{2ex}

\prop{Remark on Inner Elements}{If $y$ is an inner element of $X$
for $\langle R,X\rangle$ in FTP and connected, then $\langle R\mn
y,X\mn\{y\}\rangle$ is in FTP and connected.}

\vspace{2ex}

\noindent That $\langle R\mn y,X\mn\{y\}\rangle$ is connected is
clear from the following pictures concerning chains that ensure
connectedness:
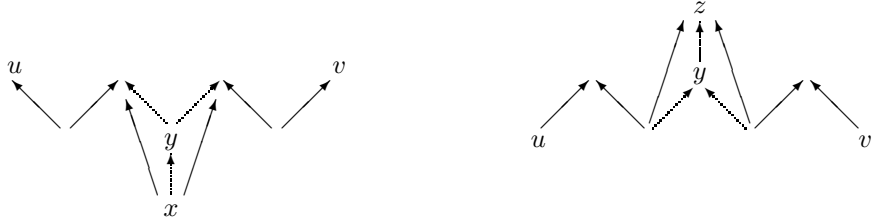
\begin{center}
\begin{picture}(330,85)(0,-15)
\put(18,20){\vector(-1,1){18}}

\put(22,20){\vector(1,1){18}}

\put(98,20){\vector(-1,1){18}}

\put(102,20){\vector(1,1){18}}

\multiput(58,22)(-1,1){14}{\line(1,0){.5}}

\multiput(62,22)(1,1){14}{\line(1,0){.5}}

\put(43,37.5){\vector(-1,1){.5}}

\put(78,37.5){\vector(1,1){.5}}

\put(-2,43){\makebox(0,0)[l]{$u$}}

\put(121,43){\makebox(0,0)[l]{$v$}}

\put(60,15){\makebox(0,0)[c]{$y$}}

\put(60,-11){\makebox(0,0)[c]{$x$}}

\multiput(59.75,-5)(0,1){12}{\line(1,0){.5}}

\put(60,10){\vector(0,1){.5}}

\put(55,-5){\vector(-1,3){12}}

\put(65,-5){\vector(1,3){12}}


\put(200,20){\vector(1,1){18}}

\put(239,20){\vector(-1,1){18}}

\put(281,20){\vector(1,1){18}}

\put(320,20){\vector(-1,1){18}}

\multiput(242,20)(1,1){14}{\line(1,0){.5}}

\multiput(277,20)(-1,1){14}{\line(1,0){.5}}

\put(262.5,35){\vector(-1,1){.5}}

\put(257,35){\vector(1,1){.5}}

\put(196,15){\makebox(0,0)[l]{$u$}}

\put(320,15){\makebox(0,0)[l]{$v$}}

\put(260,40){\makebox(0,0)[c]{$y$}}

\put(260,66){\makebox(0,0)[c]{$z$}}

\multiput(259.75,45.5)(0,1){12}{\line(1,0){.5}}

\put(260,60){\vector(0,1){.5}}

\put(279,22){\vector(-1,3){13}}

\put(241,22){\vector(1,3){13}}

\end{picture}
\end{center}
For every such chain connecting $u$ and $v$ in $\langle
R,X\rangle$ that involves the inner element $y$ there is a
substitute chain connecting $u$ and $v$ in $\langle R\mn
y,X\mn\{y\}\rangle$, which does not involve $y$. Then we have the
following.

\vspace{2ex}

\prop{Proposition 3.2}{If $\langle R,X\rangle$ is in FTP and
connected, and there are at least two elements in $X$, then for
some relations $\langle R_1,X_1\rangle$ and $\langle
R_2,X_2\rangle$ with $X_1$ and $X_2$ nonempty $\langle R,X\rangle
=\langle R_1,X_1\rangle\cdot\langle R_2,X_2\rangle$.}

\dkz We proceed by induction on the number $k$ of inner elements
of $X$ for $\langle R,X\rangle$. If $k=0$, let
\begin{tabbing}
\hspace{3em}\=$X_1 $\=$\;=\{x\in X\mid(\exists y\in X)\;(x,y)\in
R\}$,\\*[.5ex]

\>$X_2$\>$\;=\{x\in X\mid(\exists y\in X)\;(y,x)\in R\}$.
\end{tabbing}
Then $X=X_1\cup X_2$, since $\langle R,X\rangle$ is transitive and
connected, and there are at least two elements in $X$. We have
$X_1\cap X_2=\pr$, since there are no inner elements in $X$. We
also have $X_1\neq \pr$ and $X_2\neq \pr$, since $\langle
R,X\rangle$ is connected and there are at least two elements in
$X$. We can conclude that
\begin{tabbing}
\hspace{3em}\=$(\in)$\hspace{1em}\=$(\forall x_1\in X_1)(\forall
x_2\in X_2)\;(x_1,x_2)\in R$,
\end{tabbing}
because trifunctionality here implies difunctionality,\footnote{As
a matter of fact, a relation $\langle R,X\rangle$ in FTP is
difunctional iff there are no inner elements in $X$ for $\langle
R,X\rangle$.} and by difunctionality we may, to put it roughly,
shorten chains that ensure connectedness. This is clear from the
following picture:
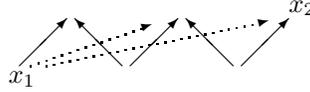
\begin{center}
\begin{picture}(110,20)(0,20)
\put(0,20){\vector(1,1){18}}

\put(39,20){\vector(-1,1){18}}

\put(42,20){\vector(1,1){18}}

\put(81,20){\vector(-1,1){18}}

\put(84,20){\vector(1,1){18}}

\put(-4,15){\makebox(0,0)[l]{$x_1$}}

\put(102,42){\makebox(0,0)[l]{$x_2$}}

\multiput(4,20)(3,1){15}{\line(1,0){.5}}

\multiput(10,20)(3,0.60){28}{\line(1,0){.5}}

\put(49.75,35.5){\vector(3,1){.5}}

\put(94,37){\vector(4,1){.5}}
\end{picture}
\end{center}
We can conclude also that
\begin{tabbing}
\hspace{3em}\=$(\notin)$\hspace{1em}\=$(\forall x_1\in
X_1)(\forall x_2\in X_2)\;(x_2,x_1)\notin R$.
\end{tabbing}
Otherwise, $\langle R,X\rangle$ would not be irreflexive. It
remains only to apply Remark~$\cdot$ of the preceding section to
establish the basis of our induction.

Suppose the number $k$ of inner elements of $X$ for $\langle
R,X\rangle$ is greater than 0. If $x$ is such an element, then, by
the Remark on Inner Elements, we have for $\langle R\mn
x,X\mn\{x\}\rangle$ too that it is in FTP and connected, and
$X-\{x\}$ has $k\mn 1$ inner elements for $\langle R\mn
x,X\mn\{x\}\rangle$. So, by the induction hypothesis, there are
relations $\langle R'_1,X'_1\rangle$ and $\langle
R'_2,X'_2\rangle$ with $X'_1$ and $X'_2$ nonempty such that
$\langle R\mn x,X\mn\{x\}\rangle =\langle
R'_1,X'_1\rangle\cdot\langle R'_2,X'_2\rangle$. We may assume that
$\langle R'_2,X'_2\rangle$ is prime with respect to $\cdot\;$, in
the sense that there are no relations $\langle S,Y\rangle$ and
$\langle T,Z\rangle$ with $Y$ and $Z$ nonempty such that $\langle
R'_2,X'_2\rangle=\langle S,Y\rangle\cdot\langle T,Z\rangle$. (If
there were such relations, then we would pass to $\langle
T,Z\rangle$ instead of $\langle R'_2,X'_2\rangle$, and rely on the
associativity of $\cdot\;$; we may iterate that.)

Since $x$ is an inner element of $X$ for $\langle R,X\rangle$,
there is a $w$ in $X$ such that $(x,w)\in R$.

1) If $w\in X'_1$, then we take
\begin{tabbing}
\hspace{3em}\=$(x1)$\hspace{1em}\=$X_1=X'_1\cup\{x\}$,\quad
$X_2=X'_2$,
\end{tabbing}
and we can conclude that $(\in)$ and $(\notin)$ hold. It remains
to apply Remark~$\cdot\;$.

2) If $w\in X'_2$, then we have the following subcases.

2.1) For every $y$ in $X'_1$ we have $(y,x)\in R$. Then we take
\begin{tabbing}
\hspace{3em}\=$(x2)$\hspace{1em}\=$X_1=X'_1$,\quad
$X_2=X'_2\cup\{x\}$,
\end{tabbing}
and we can conclude that $(\in)$ and $(\notin)$ hold. It remains
to apply Remark~$\cdot\;$.

2.2) For some element $y$ in $X'_1$ we have $(y,x)\notin R$. Then
we take $(x1)$, and to conclude that $(\in)$ and $(\notin)$ hold
it is enough to establish that
\begin{tabbing}
\hspace{3em}\=$(\ast)$\hspace{1em}\=$(\forall x_2\in
X'_2)\;(x,x_2)\in R$.
\end{tabbing}

Suppose we do not have $(\ast)$; i.e., for some $v$ in $X'_2$ we
have $(x,v)\notin R$. Let
\begin{tabbing}
\hspace{3em}\=$Y$\=$\;=\{x_2\in X'_2\mid(x,x_2)\notin
R\}$,\\*[.5ex]

\>$Z$\>$\;=\{x_2\in X'_2\mid (x,x_2)\in R\}$.
\end{tabbing}
The sets $Y$ and $Z$ are not empty, since $v\in Y$ and $w\in Z$.
Take an arbitrary $u$ from $Y$ and an arbitrary $z$ from $Z$. We
have that $(x,z)\in R$ by the definition of $Z$, and $(y,z)$ and
$(y,u)$ are in $R$ because $y\in X'_1$ and $z,u\in X'_2$. We have
$(y,x)\notin R$ by assumption, and $(x,u)\notin R$ by the
definition of $Y$. So, by trifunctionality, we may conclude that
$(u,z)\in R$, which implies $(u,z)\in R\mn x$. This implies that
for every $u$ in $Y$ and every $z$ in $Z$ we have $(u,z)\in R'_2$,
which, by Remark~$\cdot\;$, contradicts the assumption that
$\langle R'_2,X'_2\rangle$ is prime with respect to~$\cdot\;$.

So we have $(\ast)$, and hence $(\in)$ and $(\notin)$ hold. It
remains to apply Remark~$\cdot\;$. \mbox{\hspace{1em}}\qed

Then we can prove the following.

\vspace{2ex}

\prop{Proposition 3.3}{The map $K$ is onto.}

\dkz We want to show that for $\langle R,X\rangle$ in FTP there is
a diversified \textbf{S}-term $t$ such that $\kappa(t)=\langle
R,X\rangle$. We proceed by induction on the number of
\textbf{S}-variables in $X$. For the basis, if $X=\{x\}$, then
$R=\pr$, and $t$ is $x$. Suppose for the induction step that there
are at least two \textbf{S}-variables in~$X$.

If $\langle R,X\rangle$ is not connected, then, by Remark~$+$ of
the preceding section, for some relations $\langle R_1,X_1\rangle$
and $\langle R_2,X_2\rangle$ with $X_1$ and $X_2$ nonempty
$\langle R,X\rangle = \langle R_1,X_1\rangle + \langle
R_2,X_2\rangle$. So the cardinality of $X_1$ and $X_2$ is strictly
smaller than the cardinality of $X$. By Propositions 2.1 and 2.2,
we can conclude that $\langle R_1,X_1\rangle$ and $\langle
R_2,X_2\rangle$ are in FTP, and then we apply the induction
hypothesis.

If $\langle R,X\rangle$ is connected, then we apply Proposition
3.2, and reason as in the preceding paragraph. \qed

So, by the definition of $K$ and by Propositions 3.1 and 3.3, we
can conclude that $K$ is an isomorphism between a substructure of
\textbf{S} made of diversified elements and a structure on FTP.
This is an isomorphism of two algebras with partial operations $+$
and~$\cdot\;$.

\section{Shuffle sums and concatenation products on relationships}
In this and in the next section we obtain the main results of the
paper, which are summarized in the Introduction (see Section~1).
In this section we consider shuffles of arbitrary binary
relations, and with their help we define two partial operations on
sets of relations with the same domain. These operations, which we
call shuffle sum and concatenation product, are partial because we
require again disjointness of domains. The one-one map $L$, which
assigns to a partial order all its linear extensions, maps
disjoint union and concatenation of partial orders into shuffle
sum and concatenation product. With $L$, and with two other
related one-one maps, which are more general, we obtain other
isomorphic representations of the partial algebras of Section~3.

While a relation on $X$ is an ordered pair ${\langle R,X\rangle}$
such that $R\subseteq X^2$, i.e., $R\in {\cal P}(X^2)$, let a
\emph{relationship} on $X$ be an ordered pair $[U,X]$ such that
$U\subseteq {\cal P}(X^2)$, i.e., $U\in {\cal P}({\cal P}(X^2))$.
In a relationship $[U,X]$ the set $U$ is a family of the form
$\{R_i\mid i\in I\: \&\: R_i\subseteq X^2\}$. The set $X$ is the
\emph{domain} of $[U,X]$.

For the relationships $[U,X]$ and $[V,Y]$ such that $X\cap Y=\pr$
we have
\begin{tabbing}
\hspace{3em}\=$[U,X]+[V,Y]$ \= $=_{df}[\{Q\subseteq(X\cup
Y)^2\mid(\exists R\in U)(\exists S\in V)$
\\*
\` $(Q\cap X^2=R\;\;\&\;\; Q\cap Y^2=S)\},X\cup Y]$,
\\[2ex]
\>$[U,X]\,\cdot\,[V,Y]$ \> $=_{df}[\{Q\subseteq(X\cup
Y)^2\mid(\exists R\in U)(\exists S\in V)$
\\*
\` $\langle Q,X\cup Y\rangle=\langle R,X\rangle\cdot\langle
S,Y\rangle\},X\cup Y]$,
\end{tabbing}
where $\cdot$ in $\langle R,X\rangle\cdot\langle S,Y\rangle$ is
the concatenation introduced in Section~2. We call $\cdot$ in
${[U,X]\cdot[V,Y]}$, which we have just defined,
\emph{concatenation product}.

When for ${\langle R,X\rangle}$, ${\langle S, Y\rangle}$ and
${\langle Q,X\cup Y\rangle}$ such that $X\cap Y=\pr$ we have
$Q\cap X^2=R$ and $Q\cap Y^2=S$ we say that ${\langle Q,X\cup
Y\rangle}$ is a \emph{shuffle} of ${\langle R,X\rangle}$ and
${\langle S, Y\rangle}$, because this is what it is when ${\langle
R,X\rangle}$, ${\langle S, Y\rangle}$ and ${\langle Q,X\cup
Y\rangle}$ are finite linear orders. We call $+$ in $[U,X]+[V,Y]$,
defined above, \emph{shuffle sum}.

The disjoint union $\langle R,X\rangle+\langle S, Y\rangle$ and
the concatenation $\langle R,X\rangle\cdot\langle S, Y\rangle$ of
${\langle R,X\rangle}$ and ${\langle S, Y\rangle}$ are shuffles of
${\langle R,X\rangle}$ and ${\langle S, Y\rangle}$; they are limit
cases of shuffles. The disjoint union is a shuffle ${\langle
Q,X\cup Y\rangle}$ such that for every $x$ in $X$ and every $y$ in
$Y$ we have $(x,y)\notin Q$ and $(y,x)\notin Q$, while the
concatenation is a shuffle ${\langle Q,X\cup Y\rangle}$ such that
for every $x$ in $X$ and every $y$ in $Y$ we have $(x,y)\in Q$ and
$(y,x)\notin Q$ (see the Remarks $+$ and $\cdot$ in Section~2).

Consider the map $E$ from the set of relations on $X$ to the set
of relationships on $X$ defined by:
\begin{tabbing}
\hspace{3em}\=$E\langle R,X\rangle=_{df}[\{R'\subseteq X^2\mid
R\subseteq R'\},X]$.
\end{tabbing}
We use $E\langle R,X\rangle$ as an abbreviation for $E(\langle
R,X\rangle)$, and omit parentheses in the same way in analogous
situations below.

It is trivial to show that $E$ is one-one, because
$\bigcap\{R'\subseteq X^2\mid R\subseteq R'\}=R$. We can also show
that the image by $E$ of disjoint union is shuffle sum; namely, we
have the following.

\vspace{2ex}

\prop{Proposition 4.1}{$E(\langle R,X\rangle+ \langle
S,Y\rangle)=E\langle R,X\rangle+ E\langle S,Y\rangle$.}

\dkz We have to prove that $R\cup S\subseteq Q\subseteq (X\cup
Y)^2$ iff
\begin{tabbing}
\hspace{2em}\=$\exists R'\exists S'(R\subseteq R'\subseteq X^2
\;\;\&\;\; S\subseteq S'\subseteq Y^2 \;\;\&\;\; Q\cap X^2=R'
\;\;\&\;\; Q\cap Y^2=S')$.
\end{tabbing}
From left to right, it is enough to remark that from the left-hand
side we can infer that $R\subseteq Q\cap X^2\subseteq X^2$ and
$S\subseteq Q\cap Y^2\subseteq Y^2$. From right to left the
inference is trivial. \qed

On the other hand, we cannot show that $E(\langle R,X\rangle\cdot
\langle S,Y\rangle)$ is the concatenation product $E\langle
R,X\rangle\cdot E\langle S,Y\rangle$. This is because
\begin{tabbing}
\hspace{2em}\=${(Q1)}$\hspace{1em}\=$R\cup S\cup(X\times
Y)\subseteq Q\subseteq (X\cup Y)^2$
\\[1ex]
need not imply
\\[1ex]
\>${(Q2)}$\>$\exists R'\exists S'(R\subseteq R'\subseteq X^2
\;\;\&\;\; S\subseteq S'\subseteq Y^2 \;\;\&\;\; Q=R'\cup
S'\cup(X\times Y))$,
\end{tabbing}
though it is implied by it. There are sets $Q$ that satisfy
${(Q1)}$ and have in them a pair $(y,x)$ for some $x\in X$ and
some $y\in Y$.

Consider the map $P$ from the set of partial orders on $X$ to the
set of relationships on $X$ defined by replacing $R'\subseteq X^2$
in the definition of ${E\langle R,X\rangle}$ by $R'\subseteq X^2$
\emph{and} $R'$ \emph{is a partial order}. It is again trivial to
show that $P$ is one-one (for the same reason why $E$ is one-one).

Let the definitions of shuffle sum $+$ and concatenation product
$\cdot$ on relationships be modified by replacing
$Q\subseteq(X\cup Y)^2$ by $Q\subseteq(X\cup Y)^2$ \emph{and} $Q$
\emph{is a partial order}. A shuffle of two partial orders need
not be a partial order, but the concatenation of two partial
orders is a partial order (see Proposition 2.1); so the modified
definition of concatenation product amounts to the old definition
for relationships $[U,X]$ such that $U$ is a set of partial orders
on $X$. We can prove the following.

\vspace{2ex}

\prop{Proposition 4.2}{$P(\langle R,X\rangle+\langle
S,Y\rangle)=P\langle R,X\rangle+ P\langle S,Y\rangle$.}
\vspace{2ex}

\noindent For that we proceed as for Proposition 4.1. Now however
we also have the following.

\vspace{2ex}

\prop{Proposition 4.3}{$P(\langle R,X\rangle\cdot\langle
S,Y\rangle)=P\langle R,X\rangle\cdot P\langle S,Y\rangle$.}

\dkz It is enough to prove that for partial orders $Q$ the
condition ${(Q1)}$ implies ${(Q2)}$ (the converse is trivial).
Suppose ${(Q1)}$, and let $R'=Q\cap X^2$ and $S'=Q\cap Y^2$. To
show ${(Q2)}$ it is enough to show
\begin{tabbing}
\hspace{3em}\=$Q=(Q\cap X^2)\cup(Q\cap Y^2)\cup(X\times Y)$.
\end{tabbing}
To show that the right-hand side of this equation is indeed a
subset of $Q$ follows easily from ${(Q1)}$. For the converse
inclusion it is enough to verify that for every $x$ in $X$ and
every $y$ in $Y$ we cannot have $(y,x)$ in $Q$. This follows from
$X\times Y\subseteq Q$ together with the transitivity and
irreflexivity of $Q$. \qed

A relation ${\langle R,X\rangle}$ is a \emph{linear order} when it
is a partial order (as in Section~2) and for every distinct $x$
and $y$ in $X$ either $(x,y)\in R$ or $(y,x)\in R$. Consider now
the map $L$ from the set of partial orders on $X$ to the set of
relationships on $X$ defined by replacing $R'\subseteq X^2$ in the
definition of ${E\langle R,X\rangle}$ by $R'\subseteq X^2$
\emph{and} $R'$ \emph{is a linear order}. To prove that $L$ is
one-one is now not so trivial, and we need some preparation for
that.

\vspace{2ex}

\prop{Proposition 4.4}{For a partial order ${\langle R,X\rangle}$
such that for some distinct $x$ and $y$ in $X$ we have
$(y,x)\notin R$, the transitive closure ${\langle
Tr(R\cup\{(x,y)\}),X\rangle}$ is a partial order.}

\dkz We show that this transitive closure is irreflexive. If for
some $z$ in $X$ we had $(z,z)\in Tr(R\cup\{(x,y)\})$, then there
would be a chain $u_1,\ldots,u_n$ such that $u_1=u_n=z$, and
either $(u_i,u_{i+1})\in R$ or $(u_i,u_{i+1})=(x,y)$. For some $i$
we must have $(u_i,u_{i+1})=(x,y)$; otherwise $R$ would not be
irreflexive. Let $u_k$ be the leftmost $x$ in the chain, and let
$u_l$ be the rightmost $y$ in the chain. Then we must have
$(u_l,u_k)\in R$, which contradicts $(y,x)\notin R$. \qed

One can show by elementary means that every finite partial order
on $X$ can be extended to a liner order on $X$. (This is related
to what is called \emph{topological sorting} in algorithmic graph
theory.) With less elementary means one can show the same thing
for any partial order, not necessarily finite (see \cite{Je73},
p.\ 19). So, by combining this with Proposition 4.4, we obtain the
following.

\vspace{2ex}

\prop{Proposition 4.5}{For a partial order ${\langle R,X\rangle}$
such that for some distinct $x$ and $y$ in $X$ we have
$(y,x)\notin R$, there is a linear order ${\langle R',X\rangle}$
such that $R\subseteq R'$ and $(x,y)\in R'$.}

\vspace{2ex}

We can now prove that $L$ is one-one, which amounts to the
following.

\vspace{2ex}

\prop{Proposition 4.6}{For the partial orders ${\langle
R,X\rangle}$ and ${\langle S,X\rangle}$ we have that $L\langle
R,X\rangle=L\langle S,X\rangle$ implies ${R=S}$.}

\dkz Suppose $L\langle R,X\rangle=L\langle S,X\rangle$ and suppose
$(u,v)\in R$. We infer that for every linear order $S'\subseteq
X^2$ such that $S\subseteq S'$ we have $(u,v)\in S'$. If
$(u,v)\notin S$, then we obtain a contradiction with the help of
Proposition 4.5. \qed

Let the definitions of shuffle sum $+$ and concatenation product
$\cdot$ on relationships be now modified by replacing
$Q\subseteq(X\cup Y)^2$ by $Q\subseteq(X\cup Y)^2$ \emph{and} $Q$
\emph{is a linear order}. A shuffle of two linear orders need not
be a linear order, but the concatenation of two linear orders is a
linear order, and so the definition of concatenation product just
modified amounts to the old definition for relationships $[U,X]$
such that $U$ is a set of linear orders on $X$.

We can now prove the following by proceeding as for Propositions
4.1 and 4.3.

\vspace{2ex}

\prop{Proposition 4.7}{$L(\langle R,X\rangle+\langle
S,Y\rangle)=L\langle R,X\rangle+ L\langle S,Y\rangle$.}

\vspace{2ex}

\prop{Proposition 4.8}{$L(\langle R,X\rangle\cdot\langle
S,Y\rangle)=L\langle R,X\rangle\cdot L\langle S,Y\rangle$.}

\vspace{2ex}

By combining Proposition 3.1 with the facts that the maps $E$, $P$
and $L$ are one-one, we obtain new isomorphic representations of
the structure made of the diversified elements of \textbf{S} (see
Section~3).

\section{\textbf{S}-forests of graphs}
In this section we deal with the matters concerning the
constructing of graphs, which we summarized in the Introduction
(see Section~1). This is the main and concluding section of our
paper. We define first tree-like elements of the structure
\textbf{S} of Section~3, and we show that what corresponds to
these elements by the isomorphism $K$ are indeed tree-like
relations in FTP.

Consider the set $C$ of elements of \textbf{S} (see Section~3)
defined inductively as follows:
\begin{tabbing}
\hspace{1.1em}\= for every \textbf{S}-variable $x$, we have
$[x]\in C$;
\\[.5ex]
\> if $[t],[s]\in C$, then $[t+ s]\in C$;
\\[.5ex]
\> if $[t]\in C$ and $x$ is an \textbf{S}-variable, then $[x\cdot
t]\in C$.
\end{tabbing}
An alternative definition of $C$ is obtained by replacing the
third clause with:
\begin{tabbing}
\hspace{1.1em}\= if $[t]\in C$ and $+$ does not occur in the
\textbf{S}-term $s$, then $[s\cdot t]\in C$.
\end{tabbing}

Let an \textbf{S}\emph{-forest} be a diversified element of $C$.
An \textbf{S}-forest is, for example, $[((x\cdot y)\cdot z)+u]$.
An \textbf{S}\emph{-tree} is an \textbf{S}-forest that is not of
the form $[t+ s]$; for example, $[w\cdot(((x\cdot y)\cdot
z)+(u+v))]$. Since $+$ and $\cdot$ are associative, there are in
these examples superfluous parentheses, which we omit later. Note
that $[x]=\{x\}$, and that every member of $[t_1+ t_2]$ is of the
form $t_1'+ t_2'$, while every member of $[t_1\cdot t_2]$ is of
the form $t_1'\cdot t_2'$.

Let us call a partial order ${\langle R,X\rangle}$ for $X$ a
finite set of \textbf{S}-variables an FTP\emph{-forest} when for
every $x,y,z\in X$
\begin{tabbing}
\hspace{3em}\=$((x,z)\in R \;\;\&\;\; (y,z)\in R)\Rightarrow(x=y
\;\;\mbox{or}\;\; (x,y)\in R \;\;\mbox{or}\;\; (y,x)\in R)$.
\end{tabbing}
It is easy to see that FTP-forests are trifunctional, and hence
they are in FTP (see Section~3). We say that an FTP-forest
${\langle R,X\rangle}$ is an FTP\emph{-tree} when there is an
$x\in X$, called \emph{root}, such that for every $y\in X$
different from $x$ we have $(x,y)\in R$. The root is unique.
(Usually, our FTP-forests are called trees in set theory, and a
tree, which need not be finite, is defined as a partial order such
that for every element the set of its predecessors is
well-ordered.)

The following four propositions are about the map $K$ of
Section~3.

\vspace{2ex}

\prop{Proposition 5.1}{For $[t]$ an \textbf{S}-forest, $K[t]$ is
an FTP-forest.}

\dkz We proceed by induction on the length of $t$. If $t$ is an
\textbf{S}-variable, this is trivial, because $K[t]$ is the empty
relation. If $t$ is ${t_1+ t_2}$, this is trivial, by the
induction hypothesis.

Suppose $t$ is ${u\cdot t_1}$, ${(x,z)\in K[u\cdot t_1]}$ and
${(y,z)\in K[u\cdot t_1]}$. Then if $x$ and $y$ are $u$, then
${x=y}$. If $x$ is $u$ and $y$ is in $t_1$, then ${(x,y)\in
K[u\cdot t_1]}$. If $y$ is $u$ and $x$ is in $t_1$, then
${(y,x)\in K[u\cdot t_1]}$. If both $x$ and $y$ are in $t_1$, then
we apply the induction hypothesis. \qed

\prop{Proposition 5.2}{For $[t]$ an \textbf{S}-tree, ${K[t]}$ is
an FTP-tree.}

\dkz If $t$ is the \textbf{S}-variable $x$, then ${K[t]}$ is
${\langle\pr,\{x\}\rangle}$, which is an FTP-tree. If $t$ is
${[x\cdot t']}$, then $x$ is the root of ${K[t]}$. \qed

\prop{Proposition 5.3}{For every FTP-forest ${\langle R,X\rangle}$
there is an \textbf{S}-forest $[t]$ such that $K[t]=\langle
R,X\rangle$.}

\dkz By Proposition 3.3 there is a diversified \textbf{S}-term $t$
such that ${K[t]=\langle R,X\rangle}$. If $t$ has a subterm of the
form ${(s+ r)\cdot w}$, then for an \textbf{S}-variable $x$ in
$s$, an \textbf{S}-variable $y$ in $r$ and an \textbf{S}-variable
$z$ in $w$, we have ${(x,z)\in R}$, ${(y,z)\in R}$, but neither
${x=y}$, nor ${(x,y)\in R}$, nor ${(y,x)\in R}$. So ${\langle
R,X\rangle}$ is not an FTP-forest. \mbox{\hspace{1em}}\qed

\prop{Proposition 5.4}{For every FTP-tree ${\langle R,X\rangle}$
there is an \textbf{S}-tree $[t]$ such that $K[t]=\langle
R,X\rangle$.}

\dkz Just note that $t$ of the preceding proof cannot be of the
form ${t_1+ t_2}$. Otherwise ${\langle R,X\rangle}$ would not be
an FTP-tree. \qed

\noindent So $K$ establishes an isomorphism between
\textbf{S}-forests and FTP-forests on the one hand, and
\textbf{S}-trees and FTP-trees on the other hand.

We pass now to graphs and their constructing. After the following
definitions, we will give a series of examples.

A \emph{graph} is a symmetric and irreflexive relation ${\langle
G,X\rangle}$ whose domain $X$ is finite and nonempty (see
\cite{Ha69}, Chapter~2). We will now define inductively a map $T$
from the set of graphs ${\langle G,X\rangle}$ such that $X$ is a
set of \textbf{S}-variables to the power set of the set of
\textbf{S}-forests; i.e.\ ${T\langle G,X\rangle}$, which
abbreviates ${T(\langle G,X\rangle})$, is a set of
\textbf{S}-forests:

\vspace{2ex}

if ${X=\{x\}}$, then ${T\langle G,X\rangle=\{[x]\}}$;

\vspace{2ex}

\noindent supposing for the following two clauses that there are
at least two \textbf{S}-variables in $X$:

\vspace{2ex}

if ${\langle G,X\rangle}$ is connected, then
\begin{tabbing}
\hspace{3em}\=$T\langle G,X\rangle=\{[x\cdot t]\mid x\in X
\;\;\&\;\; [t]\in T\langle G\mn x,X\mn \{x\}\rangle\}$;
\end{tabbing}

if ${\langle G,X\rangle}$ is not connected, and it is of the form
${\langle G_1,X_1\rangle + \langle G_2,X_2\rangle}$ for ${\langle
G_1,X_1\rangle}$ and ${\langle G_2,X_2\rangle}$ graphs (i.e.\ for
$X_1$ and $X_2$ nonempty), then
\begin{tabbing}
\hspace{3em}\=$T\langle G,X\rangle=\{[t_1+ t_2]\mid[t_1]\in
T\langle G_1,X_1\rangle \;\;\&\;\; [t_2]\in T\langle
G_2,X_2\rangle\}$.
\end{tabbing}

It is not difficult to prove that for every graph ${\langle
G,X\rangle}$, and $x$ and $y$ distinct elements of $X$, we have
${(x,y)\in G}$ iff for every $[t]$ in ${T\langle G,X\rangle}$ the
\textbf{S}-term $t$ has no subterm ${t_1+ t_2}$ with $x$ in one of
$t_1$ and $t_2$, and $y$ in the other. From that we infer
immediately that the map $T$ is one-one.

Note that if ${\langle G,X\rangle}$ is connected, then the
\textbf{S}-forests in ${T\langle G,X\rangle}$ are
\textbf{S}-trees. These \textbf{S}-trees are in one-to-one
correspondence with what in \cite{Dev09} (Section~2) is called
maximal $(n\mn 1)$-tubings of ${\langle G,X\rangle}$, where $n$ is
the cardinality of $X$ (the notion of tubing is introduced in
\cite{CD06}, Section~2, and modified in \cite{DF08}, Section~2).
The tubings of graphs ${\langle G,X\rangle}$ that are not
connected do not however correspond exactly to the
\textbf{S}-forests in ${T\langle G,X\rangle}$.

\vspace{2ex}

\noindent \textbf{Examples 5.1.} We give now examples of
${T\langle G,X\rangle}$ for a number of connected graphs ${\langle
G,X\rangle}$.

\vspace{2ex}

\noindent \textbf{Example 5.11.} If ${\langle G,X\rangle}$ is the
connected graph
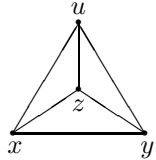
\begin{center}
\begin{picture}(50,65)
\put(0,10){\line(1,0){50}} \put(0,10){\line(3,5){25}}
\put(0,10){\line(3,2){25}} \put(50,10){\line(-3,5){25}}
\put(50,10){\line(-3,2){25}} \put(25,26.5){\line(0,1){25}}

\put(1,10){\makebox(0,0){\circle*{2}}}
\put(51,10){\makebox(0,0){\circle*{2}}}
\put(26,27){\makebox(0,0){\circle*{2}}}
\put(26,52){\makebox(0,0){\circle*{2}}}

\put(1,7){\makebox(0,0)[t]{$x$}} \put(51,7){\makebox(0,0)[t]{$y$}}
\put(25,23){\makebox(0,0)[t]{$z$}}
\put(25,55){\makebox(0,0)[b]{$u$}}
\end{picture}
\end{center}
then in ${T\langle G,X\rangle}$ we find twenty four
\textbf{S}-trees, which are obtained from the twenty four
permutations of the four \textbf{S}-variables $x$, $y$, $z$ and
$u$ by inserting $\cdot\;$. These \textbf{S}-trees naturally label
the vertices of the three-dimensional permutohedron:
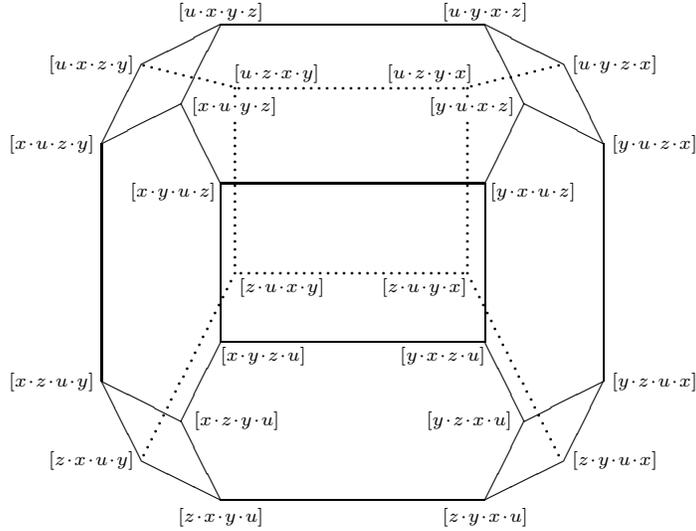
\begin{center}
\begin{picture}(200,220)(0,-10)

\put(50,10){\line(1,0){100}} \put(50,10){\line(-2,1){30}}
\put(150,10){\line(2,1){30}} \put(50,10){\line(-1,2){15}}
\put(150,10){\line(1,2){15}} \put(20,25){\line(-1,2){15}}
\put(180,25){\line(1,2){15}} \put(35,40){\line(-2,1){30}}
\put(165,40){\line(2,1){30}}

\put(35,40){\line(1,2){15}} \put(165,40){\line(-1,2){15}}
\put(50,70){\line(1,0){100}} \put(5,55){\line(0,1){90}}
\put(195,55){\line(0,1){90}} \put(50,70){\line(0,1){60}}
\put(150,70){\line(0,1){60}} \put(50,130){\line(1,0){100}}
\put(50,130){\line(-1,2){15}} \put(150,130){\line(1,2){15}}

\put(5,145){\line(1,2){15}} \put(195,145){\line(-1,2){15}}
\put(5,145){\line(2,1){30}} \put(195,145){\line(-2,1){30}}
\put(35,160){\line(1,2){15}} \put(165,160){\line(-1,2){15}}
\put(20,175){\line(2,1){30}} \put(180,175){\line(-2,1){30}}
\put(50,190){\line(1,0){100}}

\multiput(21.5,27)(1.5,3){24}{\makebox(0,0){\circle*{.5}}}
\multiput(178.5,27)(-1.5,3){24}{\makebox(0,0){\circle*{.5}}}

\multiput(21.5,174.5)(3,-.7){13}{\makebox(0,0){\circle*{.5}}}
\multiput(178.5,174.5)(-3,-.7){13}{\makebox(0,0){\circle*{.5}}}
\multiput(58.3,96)(3,0){29}{\makebox(0,0){\circle*{.5}}}
\multiput(61,166)(3,0){27}{\makebox(0,0){\circle*{.5}}}
\multiput(56,99)(0,3){19}{\makebox(0,0){\circle*{.5}}}
\multiput(56,163)(0,3){2}{\makebox(0,0){\circle*{.5}}}
\multiput(144,99)(0,3){19}{\makebox(0,0){\circle*{.5}}}
\multiput(144,163)(0,3){2}{\makebox(0,0){\circle*{.5}}}

\put(50,7){\makebox(0,0)[t]{\scriptsize $[z\!\cdot\! x\!\cdot\!
y\!\cdot\! u]$}}

\put(150,7){\makebox(0,0)[t]{\scriptsize $[z\!\cdot\! y\!\cdot\!
x\!\cdot\! u]$}}

\put(183,25){\makebox(0,0)[l]{\scriptsize $[z\!\cdot\! y\!\cdot\!
u\!\cdot\! x]$}}

\put(17,25){\makebox(0,0)[r]{\scriptsize $[z\!\cdot\! x\!\cdot\!
u\!\cdot\! y]$}}

\put(198,55){\makebox(0,0)[l]{\scriptsize $[y\!\cdot\! z\!\cdot\!
u\!\cdot\! x]$}}

\put(2,55){\makebox(0,0)[r]{\scriptsize $[x\!\cdot\! z\!\cdot\!
u\!\cdot\! y]$}}

\put(160,40){\makebox(0,0)[r]{\scriptsize $[y\!\cdot\! z\!\cdot\!
x\!\cdot\! u]$}}

\put(40,40){\makebox(0,0)[l]{\scriptsize $[x\!\cdot\! z\!\cdot\!
y\!\cdot\! u]$}}

\put(150,68){\makebox(0,0)[tr]{\scriptsize $[y\!\cdot\! x\!\cdot\!
z\!\cdot\! u]$}}

\put(50,68){\makebox(0,0)[tl]{\scriptsize $[x\!\cdot\! y\!\cdot\!
z\!\cdot\! u]$}}

\put(143,94){\makebox(0,0)[tr]{\scriptsize $[z\!\cdot\! u\!\cdot\!
y\!\cdot\! x]$}}

\put(57,94){\makebox(0,0)[tl]{\scriptsize $[z\!\cdot\! u\!\cdot\!
x\!\cdot\! y]$}}

\put(145,175){\makebox(0,0)[tr]{\scriptsize $[u\!\cdot\!
z\!\cdot\! y\!\cdot\! x]$}}

\put(55,175){\makebox(0,0)[tl]{\scriptsize $[u\!\cdot\! z\!\cdot\!
x\!\cdot\! y]$}}

\put(152,130){\makebox(0,0)[tl]{\scriptsize $[y\!\cdot\!
x\!\cdot\! u\!\cdot\! z]$}}

\put(48,130){\makebox(0,0)[tr]{\scriptsize $[x\!\cdot\! y\!\cdot\!
u\!\cdot\! z]$}}

\put(198,145){\makebox(0,0)[l]{\scriptsize $[y\!\cdot\! u\!\cdot\!
z\!\cdot\! x]$}}

\put(2,145){\makebox(0,0)[r]{\scriptsize $[x\!\cdot\! u\!\cdot\!
z\!\cdot\! y]$}}

\put(183,175){\makebox(0,0)[l]{\scriptsize $[u\!\cdot\! y\!\cdot\!
z\!\cdot\! x]$}}

\put(17,175){\makebox(0,0)[r]{\scriptsize $[u\!\cdot\! x\!\cdot\!
z\!\cdot\! y]$}}

\put(50,193){\makebox(0,0)[b]{\scriptsize $[u\!\cdot\! x\!\cdot\!
y\!\cdot\! z]$}}

\put(150,193){\makebox(0,0)[b]{\scriptsize $[u\!\cdot\! y\!\cdot\!
x\!\cdot\! z]$}}

\put(39,159){\makebox(0,0)[l]{\scriptsize $[x\!\cdot\! u\!\cdot\!
y\!\cdot\! z]$}}

\put(161,159){\makebox(0,0)[r]{\scriptsize $[y\!\cdot\! u\!\cdot\!
x\!\cdot\! z]$}}
\end{picture}
\end{center}

\noindent In this permutohedron, and in the other examples later,
there is an edge between the vertices labelled by $[t]$ and $[s]$
when there is a linear order in $L(K[t])$ and another one in
$L(K[s])$ that differ from each other just by a transposition of
immediate neighbours. (We discuss this matter after the examples.)

\textbf{Example 5.12.} If ${\langle G,X\rangle}$ is the connected
graph
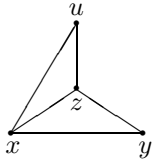
\begin{center}
\begin{picture}(50,65)
\put(0,10){\line(1,0){50}} \put(0,10){\line(3,5){25}}
\put(0,10){\line(3,2){25}} \put(50,10){\line(-3,2){25}}
\put(25,26.5){\line(0,1){25}}

\put(1,10){\makebox(0,0){\circle*{2}}}
\put(51,10){\makebox(0,0){\circle*{2}}}
\put(26,27){\makebox(0,0){\circle*{2}}}
\put(26,52){\makebox(0,0){\circle*{2}}}

\put(1,7){\makebox(0,0)[t]{$x$}} \put(51,7){\makebox(0,0)[t]{$y$}}
\put(25,23){\makebox(0,0)[t]{$z$}}
\put(25,55){\makebox(0,0)[b]{$u$}}
\end{picture}
\end{center}
obtained from the graph in the preceding example by omitting the
edge ${\{y,u\}}$, then in ${T\langle G,X\rangle}$ we find twenty
two \textbf{S}-trees, with which we label the vertices of the
following polyhedron, obtained from the three-dimensional
permutohedron by collapsing the two vertices ${[x\cdot z\cdot
y\cdot u]}$ and ${[x\cdot z\cdot u\cdot y]}$ into the single
vertex labelled ${[x\cdot z\cdot (y+ u)]}$, and the two vertices
${[z\cdot x\cdot y\cdot u]}$ and ${[z\cdot x\cdot u\cdot y]}$ into
the single vertex labelled ${[z\cdot x\cdot (y+ u)]}$:
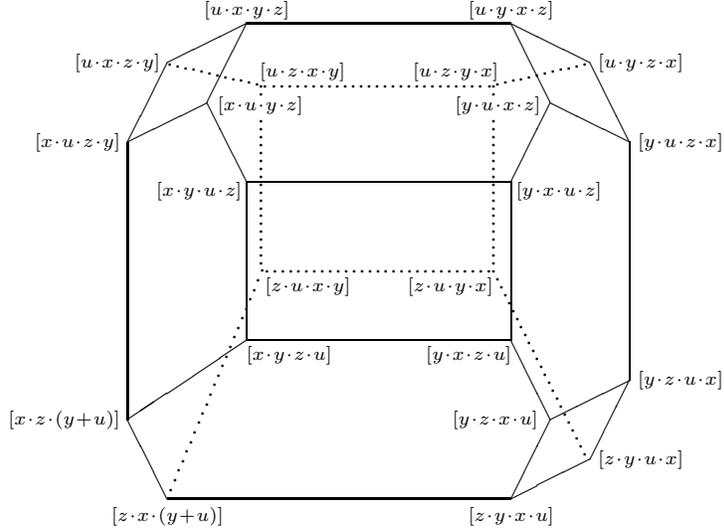
\begin{center}
\begin{picture}(200,200)

\put(20,10){\line(1,0){130}} \put(20,10){\line(-1,2){15}}
\put(50,70){\line(-3,-2){45}} \put(150,10){\line(2,1){30}}
\put(150,10){\line(1,2){15}} \put(180,25){\line(1,2){15}}
\put(165,40){\line(2,1){30}}

\put(165,40){\line(-1,2){15}} \put(50,70){\line(1,0){100}}
\put(5,40){\line(0,1){105}} \put(195,55){\line(0,1){90}}
\put(50,70){\line(0,1){60}} \put(150,70){\line(0,1){60}}
\put(50,130){\line(1,0){100}} \put(50,130){\line(-1,2){15}}
\put(150,130){\line(1,2){15}}

\put(5,145){\line(1,2){15}} \put(195,145){\line(-1,2){15}}
\put(5,145){\line(2,1){30}} \put(195,145){\line(-2,1){30}}
\put(35,160){\line(1,2){15}} \put(165,160){\line(-1,2){15}}
\put(20,175){\line(2,1){30}} \put(180,175){\line(-2,1){30}}
\put(50,190){\line(1,0){100}}

\multiput(22,13.5)(1.5,3.6){24}{\makebox(0,0){\circle*{.5}}}
\multiput(178.5,27)(-1.5,3){24}{\makebox(0,0){\circle*{.5}}}

\multiput(21.5,174.5)(3,-.7){13}{\makebox(0,0){\circle*{.5}}}
\multiput(178.5,174.5)(-3,-.7){13}{\makebox(0,0){\circle*{.5}}}
\multiput(58.3,96)(3,0){29}{\makebox(0,0){\circle*{.5}}}
\multiput(61,166)(3,0){27}{\makebox(0,0){\circle*{.5}}}
\multiput(56,99)(0,3){19}{\makebox(0,0){\circle*{.5}}}
\multiput(56,163)(0,3){2}{\makebox(0,0){\circle*{.5}}}
\multiput(144,99)(0,3){19}{\makebox(0,0){\circle*{.5}}}
\multiput(144,163)(0,3){2}{\makebox(0,0){\circle*{.5}}}

\put(20,7){\makebox(0,0)[t]{\scriptsize $[z\!\cdot\! x\!\cdot\!
(y\pl u)]$}}

\put(150,7){\makebox(0,0)[t]{\scriptsize $[z\!\cdot\! y\!\cdot\!
x\!\cdot\! u]$}}

\put(183,25){\makebox(0,0)[l]{\scriptsize $[z\!\cdot\! y\!\cdot\!
u\!\cdot\! x]$}}

\put(198,55){\makebox(0,0)[l]{\scriptsize $[y\!\cdot\! z\!\cdot\!
u\!\cdot\! x]$}}

\put(2,40){\makebox(0,0)[r]{\scriptsize $[x\!\cdot\! z\!\cdot\!
(y\pl u)]$}}

\put(160,40){\makebox(0,0)[r]{\scriptsize $[y\!\cdot\! z\!\cdot\!
x\!\cdot\! u]$}}

\put(150,68){\makebox(0,0)[tr]{\scriptsize $[y\!\cdot\! x\!\cdot\!
z\!\cdot\! u]$}}

\put(50,68){\makebox(0,0)[tl]{\scriptsize $[x\!\cdot\! y\!\cdot\!
z\!\cdot\! u]$}}

\put(143,94){\makebox(0,0)[tr]{\scriptsize $[z\!\cdot\! u\!\cdot\!
y\!\cdot\! x]$}}

\put(57,94){\makebox(0,0)[tl]{\scriptsize $[z\!\cdot\! u\!\cdot\!
x\!\cdot\! y]$}}

\put(145,175){\makebox(0,0)[tr]{\scriptsize $[u\!\cdot\!
z\!\cdot\! y\!\cdot\! x]$}}

\put(55,175){\makebox(0,0)[tl]{\scriptsize $[u\!\cdot\! z\!\cdot\!
x\!\cdot\! y]$}}

\put(152,130){\makebox(0,0)[tl]{\scriptsize $[y\!\cdot\!
x\!\cdot\! u\!\cdot\! z]$}}

\put(48,130){\makebox(0,0)[tr]{\scriptsize $[x\!\cdot\! y\!\cdot\!
u\!\cdot\! z]$}}

\put(198,145){\makebox(0,0)[l]{\scriptsize $[y\!\cdot\! u\!\cdot\!
z\!\cdot\! x]$}}

\put(2,145){\makebox(0,0)[r]{\scriptsize $[x\!\cdot\! u\!\cdot\!
z\!\cdot\! y]$}}

\put(183,175){\makebox(0,0)[l]{\scriptsize $[u\!\cdot\! y\!\cdot\!
z\!\cdot\! x]$}}

\put(17,175){\makebox(0,0)[r]{\scriptsize $[u\!\cdot\! x\!\cdot\!
z\!\cdot\! y]$}}

\put(50,193){\makebox(0,0)[b]{\scriptsize $[u\!\cdot\! x\!\cdot\!
y\!\cdot\! z]$}}

\put(150,193){\makebox(0,0)[b]{\scriptsize $[u\!\cdot\! y\!\cdot\!
x\!\cdot\! z]$}}

\put(39,159){\makebox(0,0)[l]{\scriptsize $[x\!\cdot\! u\!\cdot\!
y\!\cdot\! z]$}}

\put(161,159){\makebox(0,0)[r]{\scriptsize $[y\!\cdot\! u\!\cdot\!
x\!\cdot\! z]$}}
\end{picture}
\end{center}
We propose to call this polyhedron \emph{hemicyclohedron}. This
name will be explained in the next example. (We will not prove
here that the hemicyclohedron, conceived as an abstract polytope,
can be realized, and the same with other such polyhedra later.)

\vspace{2ex}

\noindent \textbf{Example 5.13.} If ${\langle G,X\rangle}$ is the
connected graph
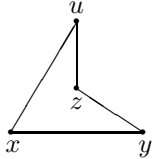
\begin{center}
\begin{picture}(50,65)
\put(0,10){\line(1,0){50}} \put(0,10){\line(3,5){25}}
\put(50,10){\line(-3,2){25}} \put(25,26.5){\line(0,1){25}}

\put(1,10){\makebox(0,0){\circle*{2}}}
\put(51,10){\makebox(0,0){\circle*{2}}}
\put(26,27){\makebox(0,0){\circle*{2}}}
\put(26,52){\makebox(0,0){\circle*{2}}}

\put(1,7){\makebox(0,0)[t]{$x$}} \put(51,7){\makebox(0,0)[t]{$y$}}
\put(25,23){\makebox(0,0)[t]{$z$}}
\put(25,55){\makebox(0,0)[b]{$u$}}
\end{picture}
\end{center}
obtained from the graph in the preceding example by omitting the
edge ${\{x,z\}}$, then in ${T\langle G,X\rangle}$ we find twenty
\textbf{S}-trees, with which we label the vertices of the
three-dimensional cyclohedron (see \cite{S97}, Section~4, and
\cite{CD06}, Corollary 2.7):
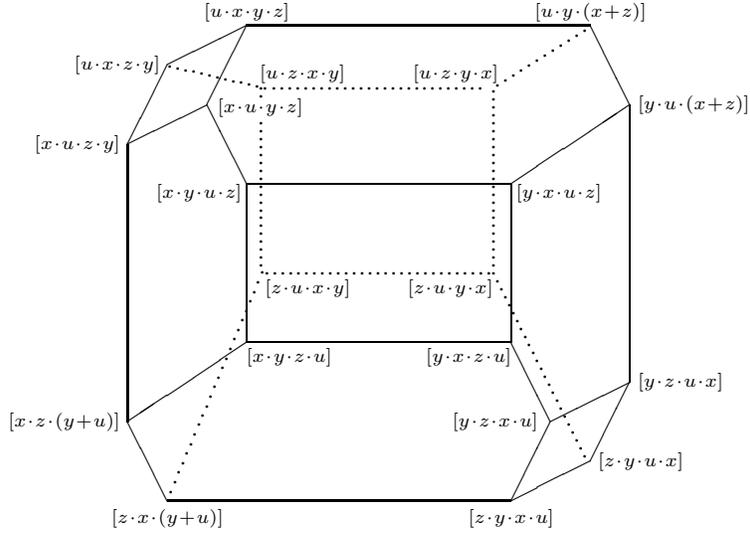
\begin{center}
\begin{picture}(200,200)

\put(20,10){\line(1,0){130}} \put(20,10){\line(-1,2){15}}
\put(50,70){\line(-3,-2){45}} \put(150,10){\line(2,1){30}}
\put(150,10){\line(1,2){15}} \put(180,25){\line(1,2){15}}
\put(165,40){\line(2,1){30}}

\put(165,40){\line(-1,2){15}} \put(50,70){\line(1,0){100}}
\put(5,40){\line(0,1){105}} \put(195,55){\line(0,1){105}}
\put(50,70){\line(0,1){60}} \put(150,70){\line(0,1){60}}
\put(50,130){\line(1,0){100}} \put(50,130){\line(-1,2){15}}
\put(150,130){\line(3,2){45}}

\put(5,145){\line(1,2){15}} \put(5,145){\line(2,1){30}}
\put(35,160){\line(1,2){15}} \put(195,160){\line(-1,2){15}}
\put(20,175){\line(2,1){30}} \put(50,190){\line(1,0){130}}

\multiput(22,13.5)(1.5,3.6){24}{\makebox(0,0){\circle*{.5}}}
\multiput(178.5,27)(-1.5,3){24}{\makebox(0,0){\circle*{.5}}}

\multiput(21.5,174.5)(3,-.7){13}{\makebox(0,0){\circle*{.5}}}
\multiput(178.5,188)(-3,-1.9){12}{\makebox(0,0){\circle*{.5}}}
\multiput(58.3,96)(3,0){29}{\makebox(0,0){\circle*{.5}}}
\multiput(61,166)(3,0){27}{\makebox(0,0){\circle*{.5}}}
\multiput(56,99)(0,3){19}{\makebox(0,0){\circle*{.5}}}
\multiput(56,163)(0,3){2}{\makebox(0,0){\circle*{.5}}}
\multiput(144,99)(0,3){23}{\makebox(0,0){\circle*{.5}}}

\put(20,7){\makebox(0,0)[t]{\scriptsize $[z\!\cdot\! x\!\cdot\!
(y\pl u)]$}}

\put(150,7){\makebox(0,0)[t]{\scriptsize $[z\!\cdot\! y\!\cdot\!
x\!\cdot\! u]$}}

\put(183,25){\makebox(0,0)[l]{\scriptsize $[z\!\cdot\! y\!\cdot\!
u\!\cdot\! x]$}}

\put(198,55){\makebox(0,0)[l]{\scriptsize $[y\!\cdot\! z\!\cdot\!
u\!\cdot\! x]$}}

\put(2,40){\makebox(0,0)[r]{\scriptsize $[x\!\cdot\! z\!\cdot\!
(y\pl u)]$}}

\put(160,40){\makebox(0,0)[r]{\scriptsize $[y\!\cdot\! z\!\cdot\!
x\!\cdot\! u]$}}

\put(150,68){\makebox(0,0)[tr]{\scriptsize $[y\!\cdot\! x\!\cdot\!
z\!\cdot\! u]$}}

\put(50,68){\makebox(0,0)[tl]{\scriptsize $[x\!\cdot\! y\!\cdot\!
z\!\cdot\! u]$}}

\put(143,94){\makebox(0,0)[tr]{\scriptsize $[z\!\cdot\! u\!\cdot\!
y\!\cdot\! x]$}}

\put(57,94){\makebox(0,0)[tl]{\scriptsize $[z\!\cdot\! u\!\cdot\!
x\!\cdot\! y]$}}

\put(145,175){\makebox(0,0)[tr]{\scriptsize $[u\!\cdot\!
z\!\cdot\! y\!\cdot\! x]$}}

\put(55,175){\makebox(0,0)[tl]{\scriptsize $[u\!\cdot\! z\!\cdot\!
x\!\cdot\! y]$}}

\put(152,130){\makebox(0,0)[tl]{\scriptsize $[y\!\cdot\!
x\!\cdot\! u\!\cdot\! z]$}}

\put(48,130){\makebox(0,0)[tr]{\scriptsize $[x\!\cdot\! y\!\cdot\!
u\!\cdot\! z]$}}

\put(198,160){\makebox(0,0)[l]{\scriptsize $[y\!\cdot\! u\!\cdot\!
(x\pl z)]$}}

\put(2,145){\makebox(0,0)[r]{\scriptsize $[x\!\cdot\! u\!\cdot\!
z\!\cdot\! y]$}}

\put(17,175){\makebox(0,0)[r]{\scriptsize $[u\!\cdot\! x\!\cdot\!
z\!\cdot\! y]$}}

\put(50,193){\makebox(0,0)[b]{\scriptsize $[u\!\cdot\! x\!\cdot\!
y\!\cdot\! z]$}}

\put(180,193){\makebox(0,0)[b]{\scriptsize $[u\!\cdot\! y\!\cdot\!
(x\pl z)]$}}

\put(39,159){\makebox(0,0)[l]{\scriptsize $[x\!\cdot\! u\!\cdot\!
y\!\cdot\! z]$}}
\end{picture}
\end{center}
Something analogous to what happened in the lower left corner of
our picture of the three-dimensional permutohedron in order to
obtain the hemicyclohedron happened now in the upper right corner
too. This explains the name of the hemicyclohedron.

\vspace{2ex}

\noindent \textbf{Example 5.14.} If ${\langle G,X\rangle}$ is the
connected graph
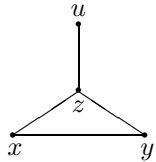
\begin{center}
\begin{picture}(50,61)(0,-1)
\put(0,10){\line(1,0){50}} \put(0,10){\line(3,2){25}}
\put(50,10){\line(-3,2){25}} \put(25,26.5){\line(0,1){25}}

\put(1,10){\makebox(0,0){\circle*{2}}}
\put(51,10){\makebox(0,0){\circle*{2}}}
\put(26,27){\makebox(0,0){\circle*{2}}}
\put(26,52){\makebox(0,0){\circle*{2}}}

\put(1,7){\makebox(0,0)[t]{$x$}} \put(51,7){\makebox(0,0)[t]{$y$}}
\put(25,23){\makebox(0,0)[t]{$z$}}
\put(25,55){\makebox(0,0)[b]{$u$}}
\end{picture}
\end{center}
obtained from the graph in Example 5.12 by omitting the edge
${\{x,u\}}$, then in ${T\langle G,X\rangle}$ we find the eighteen
\textbf{S}-trees that label the vertices of the following
polyhedron:
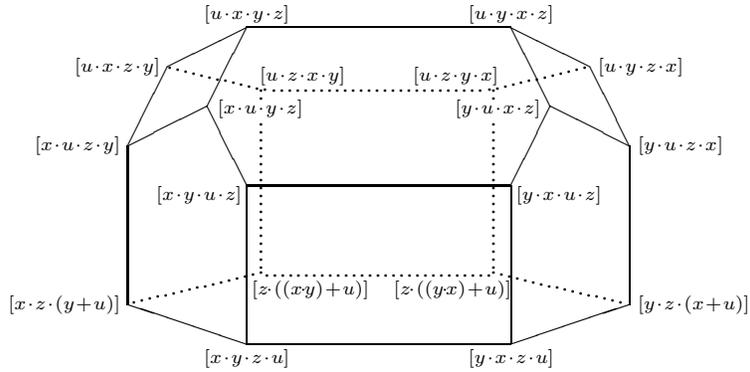
\begin{center}
\begin{picture}(200,137)(0,60)

\put(50,70){\line(-3,1){45}} \put(150,70){\line(3,1){45}}

\put(50,70){\line(1,0){100}} \put(5,85){\line(0,1){60}}
\put(195,85){\line(0,1){60}} \put(50,70){\line(0,1){60}}
\put(150,70){\line(0,1){60}} \put(50,130){\line(1,0){100}}
\put(50,130){\line(-1,2){15}} \put(150,130){\line(1,2){15}}

\put(5,145){\line(1,2){15}} \put(195,145){\line(-1,2){15}}
\put(5,145){\line(2,1){30}} \put(195,145){\line(-2,1){30}}
\put(35,160){\line(1,2){15}} \put(165,160){\line(-1,2){15}}
\put(20,175){\line(2,1){30}} \put(180,175){\line(-2,1){30}}
\put(50,190){\line(1,0){100}}

\multiput(8,85.7)(3,.7){17}{\makebox(0,0){\circle*{.5}}}
\multiput(192,85.7)(-3,.7){17}{\makebox(0,0){\circle*{.5}}}

\multiput(21.5,174.5)(3,-.7){13}{\makebox(0,0){\circle*{.5}}}
\multiput(178.5,174.5)(-3,-.7){13}{\makebox(0,0){\circle*{.5}}}
\multiput(58.3,96)(3,0){29}{\makebox(0,0){\circle*{.5}}}
\multiput(61,166)(3,0){27}{\makebox(0,0){\circle*{.5}}}
\multiput(56,99)(0,3){19}{\makebox(0,0){\circle*{.5}}}
\multiput(56,163)(0,3){2}{\makebox(0,0){\circle*{.5}}}
\multiput(144,99)(0,3){19}{\makebox(0,0){\circle*{.5}}}
\multiput(144,163)(0,3){2}{\makebox(0,0){\circle*{.5}}}

\put(198,85){\makebox(0,0)[l]{\scriptsize $[y\!\cdot\! z\!\cdot\!
(x\pl u)]$}}

\put(2,85){\makebox(0,0)[r]{\scriptsize $[x\!\cdot\! z\!\cdot\!
(y\pl u)]$}}

\put(150,68){\makebox(0,0)[t]{\scriptsize $[y\!\cdot\! x\!\cdot\!
z\!\cdot\! u]$}}

\put(50,68){\makebox(0,0)[t]{\scriptsize $[x\!\cdot\! y\!\cdot\!
z\!\cdot\! u]$}}

\put(150,94){\makebox(0,0)[tr]{\scriptsize $[z\!\!\cdot\!
((y\!\!\cdot\!\! x)\pl u)]$}}

\put(52,94){\makebox(0,0)[tl]{\scriptsize $[z\!\!\cdot\!
((x\!\!\cdot\!\! y)\pl u)]$}}

\put(145,175){\makebox(0,0)[tr]{\scriptsize $[u\!\cdot\!
z\!\cdot\! y\!\cdot\! x]$}}

\put(55,175){\makebox(0,0)[tl]{\scriptsize $[u\!\cdot\! z\!\cdot\!
x\!\cdot\! y]$}}

\put(152,130){\makebox(0,0)[tl]{\scriptsize $[y\!\cdot\!
x\!\cdot\! u\!\cdot\! z]$}}

\put(48,130){\makebox(0,0)[tr]{\scriptsize $[x\!\cdot\! y\!\cdot\!
u\!\cdot\! z]$}}

\put(198,145){\makebox(0,0)[l]{\scriptsize $[y\!\cdot\! u\!\cdot\!
z\!\cdot\! x]$}}

\put(2,145){\makebox(0,0)[r]{\scriptsize $[x\!\cdot\! u\!\cdot\!
z\!\cdot\! y]$}}

\put(183,175){\makebox(0,0)[l]{\scriptsize $[u\!\cdot\! y\!\cdot\!
z\!\cdot\! x]$}}

\put(17,175){\makebox(0,0)[r]{\scriptsize $[u\!\cdot\! x\!\cdot\!
z\!\cdot\! y]$}}

\put(50,193){\makebox(0,0)[b]{\scriptsize $[u\!\cdot\! x\!\cdot\!
y\!\cdot\! z]$}}

\put(150,193){\makebox(0,0)[b]{\scriptsize $[u\!\cdot\! y\!\cdot\!
x\!\cdot\! z]$}}

\put(39,159){\makebox(0,0)[l]{\scriptsize $[x\!\cdot\! u\!\cdot\!
y\!\cdot\! z]$}}

\put(161,159){\makebox(0,0)[r]{\scriptsize $[y\!\cdot\! u\!\cdot\!
x\!\cdot\! z]$}}
\end{picture}
\end{center}
We propose to call this polyhedron (which is called $X^a_4$ in
\cite{ACDELM09}, Figure 17) \emph{hemiassociahedron}. This name
will be explained in the next example.

\vspace{2ex}

\noindent \textbf{Example 5.15.} If ${\langle G,X\rangle}$ is the
connected graph
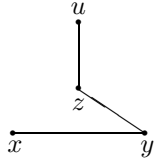
\begin{center}
\begin{picture}(50,65)
\put(0,10){\line(1,0){50}} \put(50,10){\line(-3,2){25}}
\put(25,26.5){\line(0,1){25}}

\put(1,10){\makebox(0,0){\circle*{2}}}
\put(51,10){\makebox(0,0){\circle*{2}}}
\put(26,27){\makebox(0,0){\circle*{2}}}
\put(26,52){\makebox(0,0){\circle*{2}}}

\put(1,7){\makebox(0,0)[t]{$x$}} \put(51,7){\makebox(0,0)[t]{$y$}}
\put(25,23){\makebox(0,0)[t]{$z$}}
\put(25,55){\makebox(0,0)[b]{$u$}}
\end{picture}
\end{center}
obtained from the graph in the preceding example by omitting the
edge ${\{x,z\}}$, then in ${T\langle G,X\rangle}$ we find the
fourteen \textbf{S}-trees that label the vertices of the
three-dimensional associahedron:
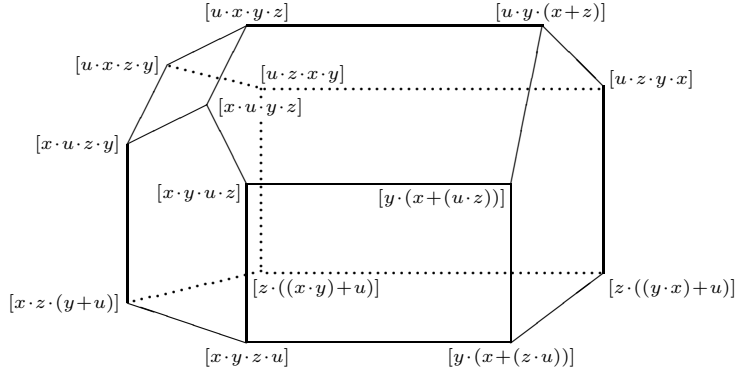
\begin{center}
\begin{picture}(200,140)(0,60)

\put(50,70){\line(-3,1){45}}

\put(50,70){\line(1,0){100}} \put(5,85){\line(0,1){60}}

\put(50,70){\line(0,1){60}} \put(150,70){\line(0,1){60}}
\put(50,130){\line(1,0){100}} \put(50,130){\line(-1,2){15}}
\put(5,145){\line(1,2){15}} \put(5,145){\line(2,1){30}}
\put(35,160){\line(1,2){15}} \put(20,175){\line(2,1){30}}
\put(50,190){\line(1,0){112.2}} \put(150,70){\line(4,3){35}}
\put(185,96){\line(0,1){71.3}}

\put(150,130){\line(1,5){11.9}} \put(162,190){\line(1,-1){23}}

\multiput(8,85.7)(3,.7){17}{\makebox(0,0){\circle*{.5}}}

\multiput(21.5,174.5)(3,-.7){13}{\makebox(0,0){\circle*{.5}}}

\multiput(58.3,96.5)(3,0){43}{\makebox(0,0){\circle*{.5}}}
\multiput(61,166)(3,0){42}{\makebox(0,0){\circle*{.5}}}
\multiput(56,99)(0,3){19}{\makebox(0,0){\circle*{.5}}}
\multiput(56,163)(0,3){2}{\makebox(0,0){\circle*{.5}}}

\put(2,85){\makebox(0,0)[r]{\scriptsize $[x\!\cdot\! z\!\cdot\!
(y\pl u)]$}}

\put(150,68){\makebox(0,0)[t]{\scriptsize $[y\!\cdot\! (x\pl
(z\!\cdot\! u))]$}}

\put(50,68){\makebox(0,0)[t]{\scriptsize $[x\!\cdot\! y\!\cdot\!
z\!\cdot\! u]$}}

\put(235,94){\makebox(0,0)[tr]{\scriptsize $[z\!\cdot\!
((y\!\cdot\! x)\pl u)]$}}

\put(52,94){\makebox(0,0)[tl]{\scriptsize $[z\!\cdot\!
((x\!\cdot\! y)\pl u)]$}}

\put(187,173){\makebox(0,0)[tl]{\scriptsize $[u\!\cdot\!
z\!\cdot\! y\!\cdot\! x]$}}

\put(55,175){\makebox(0,0)[tl]{\scriptsize $[u\!\cdot\! z\!\cdot\!
x\!\cdot\! y]$}}

\put(148,128){\makebox(0,0)[tr]{\scriptsize $[y\!\cdot\!
(x\pl(u\!\cdot\! z))]$}}

\put(48,130){\makebox(0,0)[tr]{\scriptsize $[x\!\cdot\! y\!\cdot\!
u\!\cdot\! z]$}}

\put(2,145){\makebox(0,0)[r]{\scriptsize $[x\!\cdot\! u\!\cdot\!
z\!\cdot\! y]$}}

\put(17,175){\makebox(0,0)[r]{\scriptsize $[u\!\cdot\! x\!\cdot\!
z\!\cdot\! y]$}}

\put(50,193){\makebox(0,0)[b]{\scriptsize $[u\!\cdot\! x\!\cdot\!
y\!\cdot\! z]$}}

\put(165,193){\makebox(0,0)[b]{\scriptsize $[u\!\cdot\! y\!\cdot\!
(x\pl z)]$}}

\put(39,159){\makebox(0,0)[l]{\scriptsize $[x\!\cdot\! u\!\cdot\!
y\!\cdot\! z]$}}
\end{picture}
\end{center}
In \cite{T97} it is explained how this associahedron is obtained
from the three-dimen\-sional permutohedron by two perpendicular
cuts. The previous polyhedron, the hemiassociahedron, is obtained
by one such cut. This should be also clear from our picture of the
associahedron, where one cut, which it shares with our
hemiassociahedron, is at the basis, while the other is on the
right-hand side. This explains the name of the hemiassociahedron.

\vspace{2ex}

\noindent \textbf{Example 5.16.} If ${\langle G,X\rangle}$ is the
connected graph
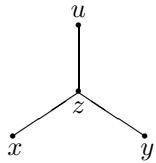
\begin{center}
\begin{picture}(50,65)
\put(0,10){\line(3,2){25}} \put(50,10){\line(-3,2){25}}
\put(25,26.5){\line(0,1){25}}

\put(1,10){\makebox(0,0){\circle*{2}}}
\put(51,10){\makebox(0,0){\circle*{2}}}
\put(26,27){\makebox(0,0){\circle*{2}}}
\put(26,52){\makebox(0,0){\circle*{2}}}

\put(1,7){\makebox(0,0)[t]{$x$}} \put(51,7){\makebox(0,0)[t]{$y$}}
\put(25,23){\makebox(0,0)[t]{$z$}}
\put(25,55){\makebox(0,0)[b]{$u$}}
\end{picture}
\end{center}
obtained from the graph in Example 5.14 by omitting the edge
${\{x,y\}}$, then in ${T\langle G,X\rangle}$ we find the sixteen
\textbf{S}-trees that label the vertices of the following
polyhedron:
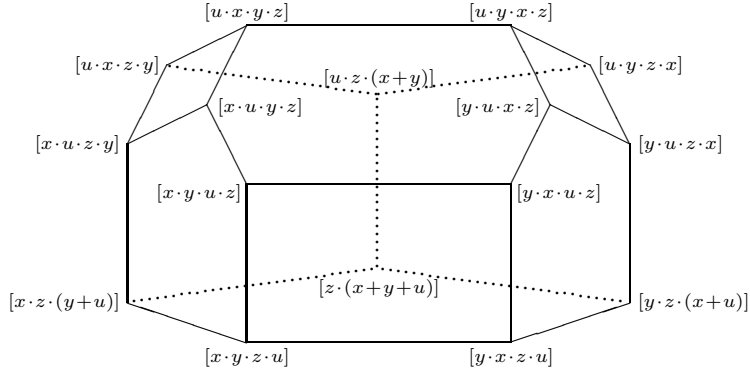
\begin{center}
\begin{picture}(200,140)(0,60)

\put(50,70){\line(-3,1){45}} \put(150,70){\line(3,1){45}}

\put(50,70){\line(1,0){100}} \put(5,85){\line(0,1){60}}
\put(195,85){\line(0,1){60}} \put(50,70){\line(0,1){60}}
\put(150,70){\line(0,1){60}} \put(50,130){\line(1,0){100}}
\put(50,130){\line(-1,2){15}} \put(150,130){\line(1,2){15}}

\put(5,145){\line(1,2){15}} \put(195,145){\line(-1,2){15}}
\put(5,145){\line(2,1){30}} \put(195,145){\line(-2,1){30}}
\put(35,160){\line(1,2){15}} \put(165,160){\line(-1,2){15}}
\put(20,175){\line(2,1){30}} \put(180,175){\line(-2,1){30}}
\put(50,190){\line(1,0){100}}

\multiput(8,85.7)(3,.4){31}{\makebox(0,0){\circle*{.5}}}
\multiput(192,85.7)(-3,.4){31}{\makebox(0,0){\circle*{.5}}}

\multiput(21.5,174.5)(3,-.4){27}{\makebox(0,0){\circle*{.5}}}
\multiput(178.5,174.5)(-3,-.4){27}{\makebox(0,0){\circle*{.5}}}

\multiput(100,99)(0,3){22}{\makebox(0,0){\circle*{.5}}}

\put(198,85){\makebox(0,0)[l]{\scriptsize $[y\!\cdot\! z\!\cdot\!
(x\pl u)]$}}

\put(2,85){\makebox(0,0)[r]{\scriptsize $[x\!\cdot\! z\!\cdot\!
(y\pl u)]$}}

\put(150,68){\makebox(0,0)[t]{\scriptsize $[y\!\cdot\! x\!\cdot\!
z\!\cdot\! u]$}}

\put(50,68){\makebox(0,0)[t]{\scriptsize $[x\!\cdot\! y\!\cdot\!
z\!\cdot\! u]$}}

\put(100,94){\makebox(0,0)[t]{\scriptsize $[z\!\cdot\! (x\pl y\pl
u)]$}}

\put(100,168){\makebox(0,0)[b]{\scriptsize $[u\!\cdot\! z\!\cdot\!
(x\pl y)]$}}

\put(152,130){\makebox(0,0)[tl]{\scriptsize $[y\!\cdot\!
x\!\cdot\! u\!\cdot\! z]$}}

\put(48,130){\makebox(0,0)[tr]{\scriptsize $[x\!\cdot\! y\!\cdot\!
u\!\cdot\! z]$}}

\put(198,145){\makebox(0,0)[l]{\scriptsize $[y\!\cdot\! u\!\cdot\!
z\!\cdot\! x]$}}

\put(2,145){\makebox(0,0)[r]{\scriptsize $[x\!\cdot\! u\!\cdot\!
z\!\cdot\! y]$}}

\put(183,175){\makebox(0,0)[l]{\scriptsize $[u\!\cdot\! y\!\cdot\!
z\!\cdot\! x]$}}

\put(17,175){\makebox(0,0)[r]{\scriptsize $[u\!\cdot\! x\!\cdot\!
z\!\cdot\! y]$}}

\put(50,193){\makebox(0,0)[b]{\scriptsize $[u\!\cdot\! x\!\cdot\!
y\!\cdot\! z]$}}

\put(150,193){\makebox(0,0)[b]{\scriptsize $[u\!\cdot\! y\!\cdot\!
x\!\cdot\! z]$}}

\put(39,159){\makebox(0,0)[l]{\scriptsize $[x\!\cdot\! u\!\cdot\!
y\!\cdot\! z]$}}

\put(161,159){\makebox(0,0)[r]{\scriptsize $[y\!\cdot\! u\!\cdot\!
x\!\cdot\! z]$}}
\end{picture}
\end{center}
Since it arises from a three-pointed star, we could perhaps call
this polyhedron the three-dimensional \emph{astrohedron}. (It is
called $D4$ in \cite{ACDELM09}, Figure 17.)

\vspace{2ex}

To sum up the previous six examples, we have the following picture
(with the number of vertices in parentheses):
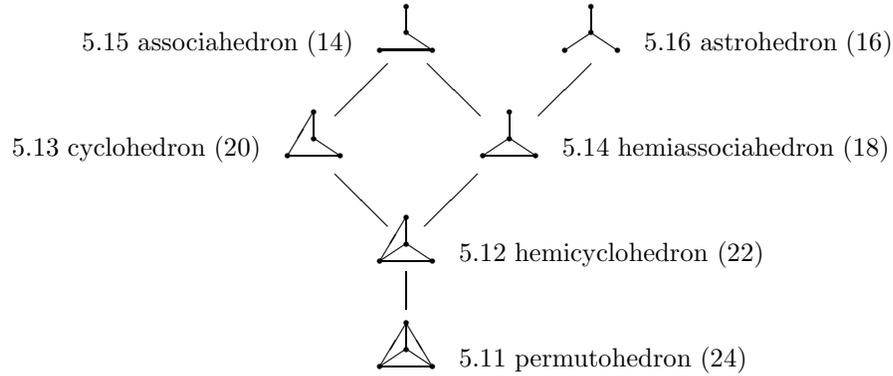
\begin{center}
\begin{picture}(120,155)
\put(40,10){\line(1,0){20}} \put(40,10){\line(3,5){10}}
\put(40,10){\line(3,2){10}} \put(60,10){\line(-3,5){10}}
\put(60,10){\line(-3,2){10}} \put(50,16.5){\line(0,1){10}}

\put(41,10){\makebox(0,0){\circle*{2}}}
\put(61,10){\makebox(0,0){\circle*{2}}}
\put(51,16.5){\makebox(0,0){\circle*{2}}}
\put(51,26.5){\makebox(0,0){\circle*{2}}}

\put(70,10){\makebox(0,0)[bl]{5.11 permutohedron (24)}}

\put(40,50){\line(1,0){20}} \put(40,50){\line(3,5){10}}
\put(40,50){\line(3,2){10}} \put(60,50){\line(-3,2){10}}
\put(50,56.5){\line(0,1){10}}

\put(41,50){\makebox(0,0){\circle*{2}}}
\put(61,50){\makebox(0,0){\circle*{2}}}
\put(51,56.5){\makebox(0,0){\circle*{2}}}
\put(51,66.5){\makebox(0,0){\circle*{2}}}

\put(70,50){\makebox(0,0)[bl]{5.12 hemicyclohedron (22)}}

\put(5,90){\line(1,0){20}} \put(5,90){\line(3,5){10}}
\put(25,90){\line(-3,2){10}} \put(15,96.5){\line(0,1){10}}

\put(6,90){\makebox(0,0){\circle*{2}}}
\put(26,90){\makebox(0,0){\circle*{2}}}
\put(16,96.5){\makebox(0,0){\circle*{2}}}
\put(16,106.5){\makebox(0,0){\circle*{2}}}

\put(-5,90){\makebox(0,0)[br]{5.13 cyclohedron (20)}}

\put(79,90){\line(1,0){20}} \put(79,90){\line(3,2){10}}
\put(99,90){\line(-3,2){10}} \put(89,96.5){\line(0,1){10}}

\put(80,90){\makebox(0,0){\circle*{2}}}
\put(100,90){\makebox(0,0){\circle*{2}}}
\put(90,96.5){\makebox(0,0){\circle*{2}}}
\put(90,106.5){\makebox(0,0){\circle*{2}}}

\put(109,90){\makebox(0,0)[bl]{5.14 hemiassociahedron (18)}}

\put(40,130){\line(1,0){20}} \put(60,130){\line(-3,2){10}}
\put(50,136.5){\line(0,1){10}}

\put(41,130){\makebox(0,0){\circle*{2}}}
\put(61,130){\makebox(0,0){\circle*{2}}}
\put(51,136.5){\makebox(0,0){\circle*{2}}}
\put(51,146.5){\makebox(0,0){\circle*{2}}}

\put(30,130){\makebox(0,0)[br]{5.15 associahedron (14)}}

\put(110,130){\line(3,2){10}} \put(130,130){\line(-3,2){10}}
\put(120,136.5){\line(0,1){10}}

\put(111,130){\makebox(0,0){\circle*{2}}}
\put(131,130){\makebox(0,0){\circle*{2}}}
\put(121,136.5){\makebox(0,0){\circle*{2}}}
\put(121,146.5){\makebox(0,0){\circle*{2}}}

\put(140,130){\makebox(0,0)[bl]{5.16 astrohedron (16)}}

{\thinlines \put(50,32){\line(0,1){14}}
\put(43,63){\line(-1,1){20}} \put(57,63){\line(1,1){20}}
\put(23,105){\line(1,1){20}} \put(78,105){\line(-1,1){20}}
\put(100,105){\line(1,1){20}} }
\end{picture}
\end{center}

We take now an example with a graph that is not connected.

\vspace{2ex}

\noindent \textbf{Example 5.2.} If ${\langle G,X\rangle}$ is the
graph
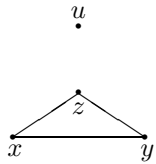
\begin{center}
\begin{picture}(50,65)
\put(0,10){\line(3,2){25}} \put(50,10){\line(-3,2){25}}
\put(0,10){\line(1,0){50}}

\put(1,10){\makebox(0,0){\circle*{2}}}
\put(51,10){\makebox(0,0){\circle*{2}}}
\put(26,27){\makebox(0,0){\circle*{2}}}
\put(26,52){\makebox(0,0){\circle*{2}}}

\put(1,7){\makebox(0,0)[t]{$x$}} \put(51,7){\makebox(0,0)[t]{$y$}}
\put(25,23){\makebox(0,0)[t]{$z$}}
\put(25,55){\makebox(0,0)[b]{$u$}}
\end{picture}
\end{center}
obtained from the graph in Example 5.14 by omitting the edge
${\{z,u\}}$, then in ${T\langle G,X\rangle}$ we find the six
\textbf{S}-forests that label the vertices of the following
hexagon:
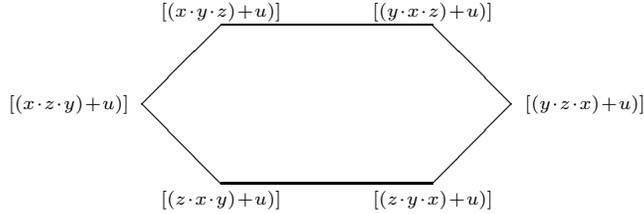
\begin{center}
\begin{picture}(200,90)

\put(60,10){\line(1,0){80}} \put(60,70){\line(1,0){80}}

\put(60,10){\line(-1,1){30}} \put(140,10){\line(1,1){30}}

\put(30,40){\line(1,1){30}} \put(170,40){\line(-1,1){30}}

\put(60,7){\makebox(0,0)[t]{\scriptsize $[(z\!\cdot\! x\!\cdot\!
y)\pl u)]$}}

\put(140,7){\makebox(0,0)[t]{\scriptsize $[(z\!\cdot\! y\!\cdot\!
x)\pl u)]$}}

\put(60,73){\makebox(0,0)[b]{\scriptsize $[(x\!\cdot\! y\!\cdot\!
z)\pl u)]$}}

\put(140,73){\makebox(0,0)[b]{\scriptsize $[(y\!\cdot\! x\!\cdot\!
z)\pl u)]$}}

\put(25,40){\makebox(0,0)[r]{\scriptsize $[(x\!\cdot\! z\!\cdot\!
y)\pl u)]$}}

\put(175,40){\makebox(0,0)[l]{\scriptsize $[(y\!\cdot\! z\!\cdot\!
x)\pl u)]$}}
\end{picture}
\end{center}

\vspace{2ex}

The \textbf{S}-forests in ${T\langle G,X\rangle}$ may be conceived
as records of the history of the destruction of ${\langle
G,X\rangle}$, which is a history of the construction of ${\langle
G,X\rangle}$ in reverse order. This destruction of graphs is based
on \emph{vertex removal} (which one finds in Ulam's Conjecture;
see \cite{Ha69}, Chapter~2). We read the \textbf{S}-forest from
left to right, and we interpret the occurrence of an
\textbf{S}-variable that we encounter in this reading as the
record of the removal of the vertex made of this
\textbf{S}-variable and of the edges involving this vertex. The
removal of vertices joined by $\cdot$ happened consecutively,
while for those joined by $+$ it happened simultaneously in time.
The commutativity of $+$ means that what is recorded on the two
sides of $+$ happened simultaneously.

For example, the \textbf{S}-forest $[(x\cdot y\cdot z)+u]$ from
Example 5.2 may be taken as a record of a destruction where,
simultaneously, one removes on the one side the vertices $x$, $y$
and $z$ and on the other side the vertex $u$; the removal of $x$,
$y$ and $z$ is done consecutively so as to produce the film:
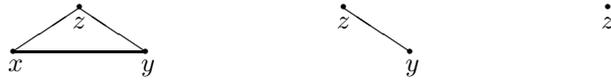
\begin{center}
\begin{picture}(250,35)
\put(0,10){\line(3,2){25}} \put(50,10){\line(-3,2){25}}
\put(0,10){\line(1,0){50}}

\put(1,10){\makebox(0,0){\circle*{2}}}
\put(51,10){\makebox(0,0){\circle*{2}}}
\put(26,27){\makebox(0,0){\circle*{2}}}

\put(1,7){\makebox(0,0)[t]{$x$}} \put(51,7){\makebox(0,0)[t]{$y$}}
\put(25,23){\makebox(0,0)[t]{$z$}}


\put(150,10){\line(-3,2){25}}

\put(151,10){\makebox(0,0){\circle*{2}}}
\put(126,27){\makebox(0,0){\circle*{2}}}

\put(151,7){\makebox(0,0)[t]{$y$}}
\put(125,23){\makebox(0,0)[t]{$z$}}


\put(226,27){\makebox(0,0){\circle*{2}}}

\put(225,23){\makebox(0,0)[t]{$z$}}
\end{picture}
\end{center}

Our examples of collapsing depend on specific graphs $\langle
G,X\rangle$, but we will show presently that we have here a
general phenomenon, not to be found only in our examples. The maps
$T$ and $L$ for a given graph $\langle G,X\rangle$ with $n$
vertices induce an equivalence relation on the set of vertices of
the $n\mn 1$-dimensional permutohedron, whose equivalence classes
are described by $L(K[t])$ for $[t]$ in $T\langle G,X\rangle$.
Moreover, the permutations corresponding to the members $L(K[t])$
are those assigned to a connected family of vertices of the
permutohedron. For example, the four permutations that correspond
to $[x\cdot y\cdot z\cdot u]$, $[x\cdot y\cdot u\cdot z]$,
$[x\cdot u\cdot y\cdot z]$ and $[u\cdot x\cdot y\cdot z]$ are
given by the linear orders in $L(K[(x\cdot y\cdot z)+u])$. The
four vertices of the permutohedron labelled by these permutations
make a connected family (see Example 5.11). We first define
precisely the required notions, and then prove three propositions,
which establish all that.

For a linear order ${\langle L,X\rangle}$ of a finite set $X$ we
call $L$ a \emph{permutation} of $X$. Let $\Lambda$ be a set of
permutations of $X$. For $L_1$ and $L_n$, where ${n\geq 2}$,
distinct members of $\Lambda$, we write ${L_1\sim_\Lambda L_n}$
when there is a sequence ${L_1\ldots L_n}$ such that
${L_1,\ldots,L_n\in\Lambda}$ and for every ${i\in\{1,\ldots,n\mn
1\}}$ we have that for two distinct $x$ and $y$ in $X$
\begin{tabbing}
\hspace{3em}\=$L_{i+1}=(L_i-\{(x,y)\})\cup\{(y,x)\}$.
\end{tabbing}
In other words, $L_{i+1}$ differs from $L_i$ just by a
transposition of immediate neighbours. We say that $\Lambda$ is
\emph{connected} when for every two distinct $L$ and $L'$ in
$\Lambda$ we have ${L\sim_\Lambda L'}$. Here are the three
propositions we announced above.

\vspace{2ex}

\prop{Proposition 5.5}{For every partial order ${\langle
R,X\rangle}$ with $X$ finite, and $L\langle
R,X\rangle=[\Lambda,X]$, the set of permutations $\Lambda$ is
connected.}

\dkz If $X$ is $\pr$ or a singleton, then ${R=\pr}$, and
${\Lambda=\{\pr\}}$, which is connected by our definition. If the
cardinality $|X|$ of $X$ is at least 2, we proceed by induction on
$|X|$.

For the basis, if ${|X|=2}$, then the only interesting case is
when $X=\{x,y\}$ and ${R=\pr}$. In that case
$\Lambda=\{\{(x,y)\},\{(y,x)\}\}$, which is clearly connected.

If ${|X|>2}$, then let $x$ be an element of $X$ such that for
every $y$ in $X$ we have ${(y,x)\notin R}$. Since $X$ is finite,
there must be such an $x$. Let $L$ and $L'$ be two different
elements of $\Lambda$. We want to show that ${L\sim_\Lambda L'}$.
Let
\begin{tabbing}
\hspace{3em}\=$S^x=\{(y,x)\mid (y,x)\in
S\}$,\hspace{5em}\=$S_x=\{(x,y)\mid (y,x)\in S\}$,
\\*[1ex]
\>$M=(L-L^x)\cup L_x$,\>$M'=(L'-L'^x)\cup L'_x$.
\end{tabbing}
It is clear that the finite sequences that correspond to the
permutations $M$ and $M'$ begin with $x$. We conclude that
${L\sim_\Lambda M}$ or ${L=M}$, and ${L'\sim_\Lambda M'}$ or
$L'=M'$. By the induction hypothesis, we have ${M\mn x\sim_\Lambda
M'\mn x}$ or ${M\mn x=M'\mn x}$. From all that ${L\sim_\Lambda
L'}$ follows. \qed

\prop{Proposition 5.6}{For every graph ${\langle G,X\rangle}$ with
$X$ a finite nonempty set of \textbf{S}-variables, and every
permutation $L$ of $X$, there is an \textbf{S}-forest $[t]$ in
${T\langle G,X\rangle}$ such that ${L\in L(K[t])}$.}

\dkz We proceed by induction on the cardinality of $X$. If $X$ is
a singleton, then we just follow the definitions. Suppose for the
induction step that $X$ has at least two \textbf{S}-variables.

If ${\langle G,X\rangle}$ is connected, then let the sequence
corresponding to the permutation $L$ be $xy_1\ldots y_n$ for
${n\geq 1}$. By the induction hypothesis, there is an
\textbf{S}-forest $s\in T\langle G\mn x,X\mn \{x\}\rangle$ such
that the permutation $L'$ corresponding to $y_1\ldots y_n$ belongs
to ${L(K[s])}$. Then we have that ${L\in L(K[x\cdot s])}$.

Suppose ${\langle G,X\rangle}$ is not connected, and is of the
form ${\langle G_1,X_1\rangle+ \langle G_2,X_2\rangle}$ for
${\langle G_1,X_1\rangle}$ and ${\langle G_2,X_2\rangle}$ graphs
(i.e.\ for $X_1$ and $X_2$ nonempty). By the induction hypothesis,
there are \textbf{S}-forests ${s_1\in T\langle G_1,X_1\rangle}$
and ${s_2\in T\langle G_2,X_2\rangle}$ such that for a permutation
$L_1$ of $X_1$ and a permutation $L_2$ of $X_2$ we have ${L_1\in
L(K[s_1])}$, ${L_2\in L(K[s_2])}$, and ${\langle L,X\rangle}$ is a
shuffle of ${\langle L_1,X_1\rangle}$ and ${\langle
L_2,X_2\rangle}$. Then we have that ${L\in L(K[s_1+ s_2])}$. \qed

\prop{Proposition 5.7}{For every graph ${\langle G,X\rangle}$ with
$X$ a finite nonempty set of \textbf{S}-variables, and every $[t]$
and $[t']$ in ${T\langle G,X\rangle}$, if $L(K[t])$ and $L(K[t'])$
are not disjoint, then ${[t]=[t']}$.}

\dkz We proceed by induction on the cardinality of $X$. If $X$ is
a singleton $\{x\}$, then ${T\langle \pr,\{x\}\rangle=\{[x]\}}$,
and ${L(K[x])=\{\langle \pr,\{x\}\rangle\}}$; the proposition
holds trivially. Suppose for the induction step that $X$ has at
least two \textbf{S}-variables.

If ${\langle G,X\rangle}$ is connected, then every element of
${T\langle G,X\rangle}$ is of the form ${[x\cdot s]}$ for some
$x$. Suppose that for some $[x\cdot s], [x'\cdot s']\in T\langle
G,X\rangle$ we have $L(K[x\cdot s])\cap L(K[x'\cdot s'])\neq\pr$.
It follows that $x$ is $x'$, and since $[s],[s']\in T\langle G\mn
x,X\mn \{x\}\rangle$ and $L(K[s])\cap L(K[s'])\neq\pr$, by the
induction hypothesis we obtain that ${[s]=[s']}$. Hence ${[x\cdot
s]=[x'\cdot s']}$.

Suppose ${\langle G,X\rangle}$ is not connected, and is of the
form ${\langle G_1,X_1\rangle + \langle G_2,X_2\rangle}$ for
${\langle G_1,X_1\rangle}$ and ${\langle G_2,X_2\rangle}$ graphs;
suppose also that $[t],[t']\in T\langle G,X\rangle$. Then, by
relying on the associativity and commutativity of $+$, we may
infer that ${[t]=[t_1+t_2]}$ and  ${[t']=[t'_1+t'_2]}$ for
$[t_i],[t'_i]\in T\langle G_i,X_i\rangle$ and ${i\in\{1,2\}}$.
Suppose for some ${\langle L,X\rangle}$ we have $\langle
L,X\rangle\in L(K[t])\cap L(K[t'])$. We infer that ${\langle L\cap
X_i^2,X_i\rangle}\in L(K[t_i])\cap L(K[t'_i])$, and hence by the
induction hypothesis we obtain that ${[t_i]=[t'_i]}$. Hence
${[t]=[t']}$. \qed

For every graph ${\langle G,X\rangle}$ with $X$ a finite nonempty
set of \textbf{S}-variables, from Propositions 5.6 and 5.7 we
infer that $\{L(K[t])\mid [t]\in T\langle G,X\rangle\}$ is a
partition of $\{\langle L,X\rangle\mid L\;\mbox{is a permutation
of}\; X\}$. We know moreover by Proposition 5.5 that every member
${L(K[t])}$ of this partition is a connected set of permutations.
Hence what we had in the examples of this section is a general
phenomenon.

\vspace{4ex}

\noindent {\small {\it Acknowledgement.} Work on this paper was
supported by the Ministry of Science of Serbia (Grants 144013 and
144029). We are grateful to Sergei Soloviev and Ralph Matthes for
useful remarks. We would like to thank also Mr.\ M.A.\ Ojlen\v
spigel for unearthing important references.}

\end{document}